\documentclass{amsart}

\usepackage{cleveref}
\usepackage{algorithm} 
\usepackage{paralist}
\usepackage{url}
\usepackage{mathrsfs}
\usepackage{color}
\usepackage{graphicx}
\usepackage{caption}
\usepackage{subcaption}
\usepackage{multirow}
\usepackage{mathtools}
\usepackage{amssymb}
\usepackage{booktabs}
\usepackage{longtable}
\usepackage{pdflscape}
\usepackage{xcolor}

\usepackage{pgfplots}
\pgfplotsset{compat=newest}
\usepackage{mathrsfs}
\usetikzlibrary{arrows}

\usepackage{tikz}
\usepackage{circuitikz}

\newcommand{\re}{\mathbb{R}}

\newcommand{\N}{\mathbb{N}}

\newcommand{\lmd}{\lambda}

\def\rank{\mbox{rank}}

\newcommand{\mc}[1]{\mathcal{#1}}

\newcommand{\st}{\mathit{s.t.}}

\newcommand{\RN}[1]{%
\textup{\uppercase\expandafter{\romannumeral#1}}%
}

\newcommand{\bdes}{\begin{description}}
\newcommand{\edes}{\end{description}}
\newcommand{\bal}{\begin{align}}
\newcommand{\eal}{\end{align}}
\newcommand{\bnum}{\begin{enumerate}}
\newcommand{\enum}{\end{enumerate}}
\newcommand{\bit}{\begin{itemize}}
\newcommand{\eit}{\end{itemize}}
\newcommand{\bea}{\begin{eqnarray}}
\newcommand{\eea}{\end{eqnarray}}
\newcommand{\be}{\begin{equation}}
\newcommand{\ee}{\end{equation}}
\newcommand{\baray}{\begin{array}}
\newcommand{\earay}{\end{array}}
\newcommand{\bsry}{\begin{subarray}}
\newcommand{\esry}{\end{subarray}}
\newcommand{\bca}{\begin{cases}}
\newcommand{\eca}{\end{cases}}
\newcommand{\bcen}{\begin{center}}
\newcommand{\ecen}{\end{center}}
\newcommand{\bbm}{\begin{bmatrix}}
\newcommand{\ebm}{\end{bmatrix}}
\newcommand{\bmx}{\begin{matrix}}
\newcommand{\emx}{\end{matrix}}
\newcommand{\bpm}{\begin{pmatrix}}
\newcommand{\epm}{\end{pmatrix}}
\newcommand{\btab}{\begin{tabular}}
\newcommand{\etab}{\end{tabular}}

\numberwithin{equation}{section}
\usepackage{lipsum}
\usepackage{amsfonts}
\usepackage{graphicx}
\usepackage{epstopdf}
\usepackage{algorithmic}
\ifpdf
\DeclareGraphicsExtensions{.eps,.pdf,.png,.jpg}
\else
\DeclareGraphicsExtensions{.eps}
\fi

\newtheorem{theorem}{Theorem}[section]

\newtheorem{prop}[theorem]{Proposition}
\newtheorem{lem}[theorem]{Lemma}
\newtheorem{lemma}[theorem]{Lemma}
\newtheorem{cor}[theorem]{Corollary}

\newtheorem{ass}[theorem]{Assumption}
\newtheorem{defi}[theorem]{Definition}
\theoremstyle{definition}
\newtheorem{example}[theorem]{Example}

\numberwithin{equation}{section}

\begin{document}

\title[PLMEs and Disjunctive Decompositions for BPOPs]
{Partial Lagrange Multiplier Expressions and Disjunctive Decompositions for Bilevel Optimization}

\author{Jiawang Nie}
\address{Jiawang Nie, Department of Mathematics,
	University of California, San Diego, 9500 Gilman Drive, La Jolla, CA, USA, 92093.}
\email{njw@math.ucsd.edu}

\author{Jane J. Ye}
\address{Jane J. Ye, Department of Mathematics and Statistics, University of Victoria,
	Victoria, B.C., Canada, V8W 2Y2.}
\email{janeye@uvic.ca}

\author{Suhan Zhong}
\address{Suhan Zhong, School of Mathematical Sciences,
	Shanghai Jiao Tong University, Shanghai, 200240, China.}
\email{suzhong@sjtu.edu.cn}

\subjclass[2020]{90C23, 65K05, 90C33}

\keywords{PLME, disjunctive decomposition, bilevel optimization,
	Lagrange multiplier, feasible extension}

\begin{abstract}
This paper studies bilevel polynomial optimization in which lower-level constraint functions
depend linearly on lower-level variables.
We show that such bilevel program can be reformulated as a disjunctive program by using
Karush-Kuhn-Tucker (KKT) conditions with a sparse type of Lagrange multipliers.
This kind of Lagrange multipliers can be conveniently represented by polynomials,
for which we call partial Lagrange multiplier expressions (PLMEs).
By doing this, each branch problem of the disjunctive program can be solved efficiently
by polynomial optimization techniques.
Solving each branch problem either returns infeasibility or gives a candidate
local or global optimizer for the original bilevel optimization.
We give necessary and sufficient conditions for these candidates to be global optimizers,
and sufficient conditions for the local optimality.
Numerical experiments are also presented to show the efficiency of the method.
\end{abstract}

\maketitle

\section{Introduction}
A typical bilevel optimization problem (BOP) is in the form
\begin{equation}\label{eq:genblv}
	\left\{
	\begin{array}{cl}
		\min\limits_{x\in\mathbb{R}^n,y\in\mathbb{R}^p} & F(x,y)\\
		\st & h_i(x,y)\ge 0\,(i\in\mathcal{I}_1),\\
			& y\in S(x),
	\end{array}
	\right.
\end{equation}
where $S(x)$ is the set of global optimizer(s) of the lower-level problem
\[
	(P_x)\quad\left\{
	\begin{array}{cl}
		\min\limits_{y\in\mathbb{R}^p} & f(x,y)\\
		\st &  g_j(x,y) \ge 0\,(j\in\mathcal{I}_2).
	\end{array}\right.
\]
Here $\mathcal{I}_1,\mathcal{I}_2$ are finite (or empty) index sets.
We remark that an equality constraint $h=0$ can be equivalently given as
two inequality constraints $h\ge 0$ and $h\le 0$.
We call problem \eqref{eq:genblv} a bilevel polynomial optimization problem (BPOP)
if all $F, f, h_i, g_j$  are polynomials.
The functions $F$ and $h_i$ are called the upper-level objective function
and constraint functions respectively,
while $f$ and $g_j$ are called the lower-level objective function
and constraint functions respectively.

Let $Z(x)$ be the feasible set of the lower-level problem $(P_x)$ and denote
\begin{equation*}
	\mathcal{U} \coloneqq
	\{ (x,y)\,\vert\, h_i(x,y)\ge 0\,(i\in \mathcal{I}_1),\,
	      g_j(x,y )\ge 0 \, (j\in\mathcal{I}_2) \}.
\end{equation*}
The feasible set of \eqref{eq:genblv} can be written as
\begin{equation*}
	\mc{F} \coloneqq \{ (x,y)\in \mc{U}\,\vert\, y\in S(x) \}.
\end{equation*}
Throughout the paper we assume that the solution set $S(x)$ is nonempty
for all feasible $x$ satisfying  $(x,y)\in \mc{U}$ for some $y$.

Bilevel optimization originated from economics as
Stackelberg games \cite{Mirrlees99,Stackelberg10}.
It has been successfully applied to many application areas,
and more recently to data science and machine learning (see e.g., \cite{Franceschi,GaoYeYinZengZhang,Kunapuli08,YeDC2021}).
We refer to the monographs \cite{Bardbook13,Dempebook18,dempe2015bilevel,Shimizu97},
the surveys \cite{ColsonSuv07,DempeZemkoho20}
and the references within for more applications and the recent advances in related topics.

To study a bilevel optimization problem,
one usually transforms it into a single-level optimization problem.
However, the equivalent single-level program is still very challenging to solve.
This is because the usual constraint qualifications always fail at every point
in the feasible region,
and implicit functions like the value function may be used in the reformulation of
bilevel optimization \cite{SBLYe2013,ye1995optimality}.
Despite these difficulties, there are still tremendous developments on
constraint qualifications and optimality conditions for bilevel optimization.
We refer interested readers to \cite{BaiYe21,DempeZemkoho20,KeYaoYeZhang,MaYaoYeZhang,yecq20}
and the references within for recent advances about these topics.

It is notoriously hard to design an efficient and reliable algorithm to solve
bilevel optimization problems.
Even when all defining functions are linear,
the computational complexity is already NP-hard \cite{AyedBlair90}.
In most early works in the literature,
bilevel optimization problems are tackled by replacing the lower-level problem
with its first-order optimality conditions like the Karush-Kuhn-Tucker (KKT) conditions.
This KKT reformulation then minimizes over the original variables and the Lagrange multipliers.
When the lower-level constraints include inequality constraints,
the transformed problem becomes a so-called mathematical program with
equilibrium constraints (MPEC) \cite{LuoPangRalph,Outrata}.
But this KKT approach is only reliably applicable if the lower-level problem is smooth and convex.
Moreover, when the lower-level multipliers are not unique,
the resulting optimization usually has different local optimizers
from the original ones \cite{Dempe12lo,LamSag20}.
In recent years, some numerical algorithms for bilevel optimization that
are not reformulated as MPECs are proposed in
\cite{GaoYeYinZengZhang,LamSag20,SBLYe2013,mitsos2008global,outrata1990numerical,xu2013smoothing,
XuYeZhang,YeDC2021}.

In this paper, we mainly focus on solving bilevel polynomial optimization
whose lower-level constraining functions are linear in lower-level variables.
Consider
\begin{equation}\label{eq:lin:g_i}
	g_j(x,y) \, = \, a_j^Ty-b_j(x),\quad \text{for} \quad j\in\mc{I}_2,
	\end{equation}
where each $a_j$ is a constant vector and $b_j(x)$ is a scalar polynomial.
This class of bilevel optimization problem is broadly applied in
portfolio design \cite{Gon2021}, price setting problems \cite{Labbe16},
toll optimization \cite{Kue2023}.
In this paper, we consider a toll-setting problem in Example~\ref{ex:toll}.
BPOPs can be relaxed into single-level deterministic polynomial optimization.
Polynomial optimization has been extensively studied
\cite{Las01,Las09book,Lau09,nie2014optimality,MPO}.
A single-level polynomial optimization problem can be solved globally
by Moment-SOS relaxations.
Recently, these techniques have been used to design numerical algorithms
to globally solve BPOPs in \cite{Jeyakumat16,nie2017bilevel,nie2020bilevel}.

\subsection{Lagrange multiplier expression approach}
\label{ssec:lme-fe}

For convenience, we assume
\[
	\mc{I}_2 = [m]\coloneqq \{1,\ldots, m\},
\]
and denote the vector of lower-level constraining functions
\[
	g(x,y) \coloneqq  Ay-b(x)\quad \mbox{where}\quad
\]
\[
	A \coloneqq \bbm a_1 &  \cdots & a_m\ebm^T,\quad
	b(x)  \coloneqq \bbm b_1(x) & \cdots &  b_m(x)\ebm^T .
\]
Every $(x,y)\in \mc{F}$ satisfies the KKT conditions:
\be \label{eq:KKT:intro}
	\exists \,  \lambda=(\lambda_1,\ldots, \lambda_m)  \quad \st\quad
	\left\{ \begin{array}{l}
	\nabla_y f(x,y) -A^T \lambda =0,\\
	0 \le [a_j^Ty-b_j(x)] \perp \lambda_j \geq 0 \ ( j\in [m]),
	\end{array}
	\right.
\ee
where $\perp$ denotes the perpendicularity and
the $\lambda_j$'s are Lagrange multipliers.
Denote the KKT set
\begin{equation*}
	\mc{K} \coloneqq \left\{ (x,y)\in\mc{U} \left|
	\exists \,\, \lambda \quad \st \quad
	\begin{array}{l}
		\nabla_y f(x,y) -  A^T\lambda= 0,\\
		0 \le [a_j^Ty-b_j(x)] \perp \lambda_j \geq 0  \ ( j\in [m])
	\end{array}
	\right.\right\}.
\end{equation*}
By replacing $y\in S(x)$ with the KKT conditions \eqref{eq:KKT:intro}
and the KKT set constraint $(x,y)\in \mc{K}$ respectively,
one can obtain relaxations of \eqref{eq:genblv} as
\[\begin{array}{l|l}
	(P_{kkt}):\,
		\left\{\begin{array}{cl}
		\min\limits_{(x,y,\lambda)} & F(x,y)\\
		\st & (x,y)\in \mc{U},\\
			& (x,y,\lambda)\, \mbox{satisfies}\, \eqref{eq:KKT:intro}
		\end{array}
	\right. &
	(P):\,
		\left\{\begin{array}{cl}
		\min & F(x,y)\\
		\st & (x,y)\in \mc{K}.
		\end{array}\right.
\end{array}\]
The $(P_{kkt})$ is a polynomial optimization problem.
Suppose $(P)$ is also defined by polynomials.
Then $(P_{kkt})$ and $(P)$ can be solved globally by Moment-SOS relaxations
(if solved by other methods,
the global optimality may not be guaranteed; see in \Cref{tab:compare}).
The dimensions of the decision variables for the $l$-th order moment relaxations of
$(P_{kkt})$ and $(P)$, which take the form of \eqref{dpr},
are respectively
\[
	\binom{n+p+m+2l}{2l}\quad
	\mbox{and}\quad
	\binom{n+p+2l}{2l}.
\]
Clearly, when $m$ is big, the moment relaxation of $(P_{kkt})$
has much bigger size than $(P)$, for the same relaxation order.
So it is more efficient to design an algorithm for BPOPs
based on $(P)$, especially when there are many lower-level constraints.

The strategy recently proposed in \cite{nie2020bilevel}
can be used to design such an algorithm.
Suppose there is a vector $\lambda(x,y)$  of rational functions
\footnote{A rational function is in form of the ratio of two polynomials.}
such that
\[
	\mbox{$(u,v,\lambda(u,v))$ satisfies \eqref{eq:KKT:intro}} \quad
	\mbox{for every}\quad (u,v)\in\mc{K}.
\]
Such a $\lambda(x,y)$ is called a {\it Lagrange multiplier expression} (LME).
Then $(P)$ can be approximated by the rational optimization
\begin{equation}\label{eq:KKTrel:LME}
	\left\{
	\begin{array}{cl}
		\min\limits_{(x,y)\in\mc{U}} & F(x,y)\\
		\st & \nabla_y f(x,y) - A^T\lambda(x,y) = 0,\\
		& 0 \le [a_j^Ty-b_j(x)]\perp \lambda_j(x,y) \geq 0 \
		( j\in \mathcal{I}_2).
	\end{array}
	\right.
\end{equation}
Suppose $(\hat{x}, \hat{y})$ is the global minimizer of problem \eqref{eq:KKTrel:LME}.
If $\hat{y}\in S(\hat{x})$, then $(\hat{x}, \hat{y})$ is also the
global optimizer of \eqref{eq:genblv}.
Otherwise, one can find $\hat{z}\in Z(\hat{x})$ such that
$f(\hat{x}, \hat{z})<f(\hat{x},\hat{y})$.
Assume there is a vector-valued rational function $q:\re^n\times \re^p\to \re^p$
that satisfies
\begin{equation}\label{eq:fe}
	q(\hat{x},\hat{y}) = \hat{z},\quad\mbox{and}\quad
	q(x,y)\in Z(x)\,\,  \forall (x,y)\in \mc{K}.
\end{equation}
Such a $q(x,y)$ is called a {\it feasible extension} (FE) of $\hat{z}$
at $(\hat{x},\hat{y})$. If we add
\begin{equation}\label{eq:diff_FE}
	f(x,y) - f(x,q(x,y))\ge 0
\end{equation}
as an extra constraint to \eqref{eq:KKTrel:LME},
then $(\hat{x},\hat{y})$ will no longer be feasible and
\eqref{eq:KKTrel:LME} will become a tighter relaxation of $(P)$
with a different global optimizer.
Under the existence of feasible extensions and some other continuity assumptions,
this process can be repeated infinitely and the computed candidate solutions
converge to the global optimizer of \eqref{eq:genblv}.
Recently, this LME-FE approach was applied to solve BPOPs \cite{nie2020bilevel}
and general Nash equilibrium problems \cite{NieTangZgnep21}.

The LME-FE approach exhibits evident computational advantages compared to
the classic KKT approach (see comparison in \cite[Example~6.1]{nie2020bilevel}).
It prevents Lagrange multipliers being extra decision variables and obtains
a strict hierarchy of outer approximations (possibly tight) for $\mc{F}$.
The LME-FE approach is appealing for finding optimizers of BPOPs,
especially for nonconvex cases.
LMEs exist for generic polynomial constraints (see \Cref{ssc:exLME})
and FEs have universal expressions for boxed and simplex constraints.
However, they usually have very complicated expressions (e.g., with high degrees)
when the lower-level problem has a relatively large number of constraints.
This makes \eqref{eq:KKTrel:LME} be computationally expensive to solve.
For general linear constraints, we are interested in:
\begin{itemize}
	\item How to automatically find more efficient LMEs?
	\item How to identify the existence and compute convenient expressions of FEs?
\end{itemize}
These issues are addressed in this paper.

\subsection{Disjunctive decomposition}

To reduce the complexity, we consider finding more convenient LMEs to decompose
the lower KKT set such that each decomposed piece is easy to represent.

Let $J = \{j_1,\ldots, j_t\}$ be a subset of $[m]$ with cardinality $t$, where
\[
	t \coloneqq \rank\, A.
\]	
Consider Lagrange multiplier $\lambda$ such that $\lambda_j = 0$ for all $j\not\in J$.
Then the KKT conditions in \eqref{eq:KKT:intro} become
\be \label{kktcond}
	\left\{ \baray{rrl}
	\exists \,  \lambda  \quad \st \quad &
	\nabla_y f(x,y) -  A^T \lambda=0, &   \\
	& \lambda_j =0\, & (j\not \in J), \\
	& 0 \le [a_j^Ty-b_j(x)] \perp \lambda_j \geq 0\,
	& ( j\in J).
	\earay \right.
\ee
Computing LMEs from \eqref{kktcond} is much easier than \eqref{eq:KKT:intro}
with priori zero entries. Let
\[
	\mc{K}_J   \, \coloneqq  \, \big \{ (x,y)\in\mc{U} \,\vert\,\,
	\exists \,\, \lambda \,\,
	\text{satisfying} \, \eqref{kktcond} \big \}.
\]
We define the $J$-th {\it branch problem} of $(P)$ as
\[
	(P_J):\quad \left\{\begin{array}{cl}
	\min & F(x,y)\\
	\st & (x,y)\in \mc{K}_J.
	\end{array}
	\right.
\]
If $\{a_j\,\vert\, j\in J\}$ is linearly independent,
then $\mc{K}_J$ can be expressed as a semialgebraic set
with {\it partial Lagrange multiplier expressions} (PLMEs)
given by polynomials (see \Cref{prop:PLME}).
As shown in \Cref{thm:kktdcp}, every KKT point of \eqref{eq:genblv} is
contained in such a $\mc{K}_J$.
So we can find the solution of $(P)$ by solving all branch problems and
selecting the optimizer that achieves the minimum value.
In particular, $(P_J)$ is usually easy to solve since PLMEs
are convenient to obtain.
When the lower-level problem is convex,
every KKT point is feasible for the original bilevel optimization,
then $(P)$ and \eqref{eq:genblv} have the same optimizer(s).
When $f$ is nonconvex, we apply the FE approach introduced in the previous discussions.

Solving an optimization by its branches uses the idea of {\it disjunctive programming}.
This approach was introduced by Balas \cite{Balas} to represent mixed-integer
programming problems. We refer to \cite{Lee,Raman} for related work.
Recently, disjunctive programming methods based on KKT reformulations
for MPECs are studied in \cite{Liang,Mehlitz}.
We would like to remark that our approach is quite different from these methods.
First, each $(P_J)$ does not have Lagrange multipliers as decision variables.
More importantly, in each branch problem,
we search for an optimizer that is also feasible for the original BPOP
using feasible extensions.
Hence our method can be used to get the local/global optimizers and
verify their optimality.

\subsection*{Contributions}

Our main contributions are the following.
\begin{itemize}

\item
We propose the disjunctive program reformulation to solve bilevel programs.
We remark that our disjunctive program is not enumerating active sets
for the lower-level optimization.
This is because some constraints for the label set $J$ may still be inactive.
With techniques of PLMEs and feasible extensions,
it also has major differences from disjunctive program from
the KKT reformulation $(P_{kkt})$.
To the best of our knowledge, this novel approach has never been
used to solve bilevel optimization before.

\item
We propose to solve BPOPs with linear lower-level constraints as in \eqref{eq:lin:g_i}
by the disjunctive decomposition approach combined with the
Lagrange multiplier expression approach.
In particular, the construction of disjunctive program does
not require the linearity of lower-level constraints.
It is a promising future work to study the disjunctive decomposition approach
for BPOPs with more general constraints.

\item
We give convenient sufficient conditions to verify the local optimality,
and sufficient and necessary conditions for verifying global optimality of BPOPs.
The solution candidates are solved from branch problems by
Moment-SOS relaxations with feasible extensions.

\end{itemize}

\subsection*{Notation}

The symbol $\mathbb{N}$  (resp., $\mathbb{R}$, $\re_+$)
denotes the set of nonnegative integers (resp., real numbers, nonnegative real numbers).
For an integer $m>0$, $[m] \coloneqq \{1,\cdots,m\}$.
For $t\in \re$, $(t)_+ \coloneqq \max\{t,0\}$ and $\lceil t\rceil$ denotes
the smallest integer that is greater than $t$.
For a function $f$, $\nabla_y f$ denotes its partial gradient with respect to $y$
and $\nabla_{y_i}f$ denotes the derivative of $f$ in $y_i$.
The notation $\| \cdot \|$ denotes the standard Euclidean norm for a vector or a matrix.
The $e \coloneqq (1,\ldots,1)^T$ denotes the vector of all ones and
each $e_i$ denotes the unit vector of all zeros but the $i$-th entry being one.
The symbol $0_{n_1\times n_2}$ denotes the zero matrix of dimension $n_1\times n_2$
and $I_m$ denotes the identity matrix of dimension $m$.
A symmetric matrix $A\in\re^{n\times n}$ is said to be positive semidefinite
or psd if $x^TAx\ge 0$ for every $x\in\re^n$, denoted as $A\succeq 0$.
The closed ball centered at $x$ with radius $\rho>0$ is denoted by $B_\rho(x)$.
Let $J\subseteq [m]$.
For a matrix $A=[a_1\,\,\cdots\,\, a_m]^T$, we denote $A_J \coloneqq [a_j^T]_{j\in J}$
as the submatrix of $A$ with rows labelled by $J$.
For a vector $b = (b_1,\ldots, b_m)^T$, we denote $b_J \coloneqq (b_j)_{j\in J}$.

The rest of paper is organized as follows.
In \Cref{sc:DecpKKT}, we introduce partial Lagrange multiplier expressions and decomposed KKT system.
In \Cref{sc:disjunctive}, we propose the disjunctive decomposition of BPOPs.
In \Cref{sec:alglocal}, we propose semidefinite algorithms to solve BPOPs.
The construction of feasible extensions, as a key component of the algorithms,
is studied in \Cref{sc:FE}.
The numerical examples are presented in \Cref{sc:numexam}.
Conclusions are made in \Cref{sc:con}.

\section{PLMEs and the KKT system}
\label{sc:DecpKKT}

Recall that the KKT set of \eqref{eq:genblv} is
\begin{equation}\label{eq:kktset}
	\mc{K} \coloneqq \left\{ (x,y)\in\mc{U} \left|\,
	\exists \,\, \lambda \quad \st \quad
	\begin{array}{l}
		\nabla_y f(x,y) -  A^T\lambda= 0,\\
		0 \le [a_j^Ty-b_j(x)] \perp \lambda_j \geq 0  \ ( j\in [m])
	\end{array}
	\right.\right\}.
\end{equation}

\subsection{Existence of LMEs}\label{ssc:exLME}

First, we discuss the existence of Lagrange multiplier expressions (LMEs).
Recall that a vector-valued rational function $\lambda(x,y)$ is a LME if
\[
	\left\{\begin{array}{l}
	\nabla_y f(x,y) - A^T\lambda(x,y) = 0,\\
	0\le [a_j^Ty-b_j(x)] \perp \lambda_j(x,y)\ge 0\,\,(j\in[m]).
	\end{array}\right.
\]
is satisfied at every $(u,v)\in\mc{K}$. Denote
\[
	G(x,y) = \bbm a_1 & a_2 & \cdots & a_m\\
	g_{1}(x,y) & 0 & \cdots & 0\\
	\vdots & \vdots & & \vdots\\
	0 & 0 & \cdots & g_m(x,y)\ebm,\quad
	\hat{f}(x,y) = \bbm \nabla_y f(x,y)\\ 0\\\vdots\\ 0\ebm.
\]
A LME $\lambda(x,y)$ is necessary to satisfy
\begin{equation}\label{eq:Glmd=f}
	G(x,y)\lambda(x,y) =\hat{f}(x,y).
\end{equation}
If there exist a matrix of polynomials $L(x,y)$ and a nonzero scalar polynomial
$d(x,y)$ such that
\begin{equation}\label{eq:LG=dI}
	L(x,y)G(x,y) = d(x,y) I_m,
\end{equation}
then $d(x,y)\lambda(x,y) = L(x,y)\hat{f}(x,y)$ by multiplying \eqref{eq:Glmd=f}
with $L(x,y)$. It implies the LME given by rational functions:
\begin{equation} \label{eq:lme}
	\lambda(x,y) = L(x,y)\hat{f}(x,y)/d(x,y),
\end{equation}
which satisfies $L(x,y)\hat{f}(x,y) = 0$ if $d(x,y) = 0$.
Such $L(x,y), d(x,y)$ exist for general cases.
For example, let $H(x,y)\coloneqq G(x,y)^TG(x,y)$. Consider
\begin{equation}\label{eq:Lnd}
	L(x,y) = \operatorname{adj} (H(x,y))G(x,y)^T,\quad
	d(x,y) = \det(H(x,y)),
\end{equation}
where $\operatorname{adj}(\cdot)$ denotes the matrix adjugate
and $\det(\cdot)$ is the determinant notation.
Since $\det (H(x,y))$ is not identically zero in general cases, we have
\[
	\operatorname{adj}(H(x,y)) H(x,y) = \det(H(x,y)) I_m.
\]
The condition \eqref{eq:LG=dI} is satisfied with polynomials given in \eqref{eq:Lnd}.
We remark that there usually exist different options for $L(x,y)$ and $d(x,y)$
to make \eqref{eq:LG=dI} hold.
We prefer those with low degrees for computational efficiency.
However, we have not found simpler universal expressions other than \eqref{eq:Lnd},
except in special cases such as boxed and simplex constraints \cite{nie2017bilevel}.

\subsection{Partial Lagrange multiplier expressions}

Recall $\rank\, A = t$.
Let $J = \{j_1,\ldots, j_t\}\subseteq [m]$. Denote that
\[
	A_J = [a_{j_1}\,\, \cdots\,\, a_{j_t}]^T,\quad
	b_J(x) = [b_{j_1}(x),\,\,\cdots\,\,b_{j_t}(x)]^T.
\]
We consider the Lagrange multiplier $\lmd$ such that $\lmd_j = 0$
for all $j \not\in J$.
It determines the KKT subset associated with $J$ as
\be \label{eq:defKJ}
	\mc{K}_J \coloneqq \left\{ (x,y)\in \mc{U}\left|\baray{l}
	\exists\,  \lambda_J = (\lambda_{j_1},\ldots,\lambda_{j_t})\quad \st\\
	\nabla_y f(x,y) -  A_J^T \lambda_J = 0,   \\
	0 \le [a_j^Ty-b_j(x)] \perp \lambda_j \geq 0\,\,( j\in J)
	\earay \right.\right\}.
\ee
We can similarly define Lagrange multiplier expressions for such a $\mc{K}_J$.

\begin{defi} \label{df:PLME}
	For $J = \{ j_1, \ldots, j_t \}\subseteq [m]$,
	a function $\psi:\re^n\times \re^p\to \re^t$ is called a
	partial Lagrange multiplier expression (PLME) with respect to $J$ if
	\begin{equation}\label{eq:def:plme}
		A_J^T \psi(x,y) =
		\nabla_y f(x,y)\quad \mbox{for all}\quad (x,y) \in \mc{K}_J.
	\end{equation}
	For such a case, we write $\lmd_J(x,y)\coloneqq \psi(x,y)$ with each
	$\lmd_{j_i}(x,y) \coloneqq \psi_{i}(x,y)$.
\end{defi}

The function $\psi$ satisfying Definition~\ref{df:PLME} always exists
if $A_J$ has full rank. We denote the collection of such $J$ as
\begin{equation}\label{set:Pmt}
	\mc{P} \coloneqq \{J\subseteq [m]\,\vert\,
	|J| = t,\, A_J\,\mbox{has full rank} \}.
\end{equation}
Let $J\in \mc{P}$. Since $t = \rank\, A\le \min\{m,p\}$,
the full rank matrix $A_J\in\re^{t\times p}$ has linearly independent rows,
thus $A_JA_J^T$ is invertible. From the KKT equation
\[	A_J^T\lambda_J = \nabla_y f(x,y),	\]
we can solve for the PLME
\begin{equation}\label{eq:PLMEgen}
	\lambda_J(x,y) = (A_JA_J^T)^{-1}A_J\nabla_y f(x,y).
\end{equation}
In particular, if $t=p$, then $A_J$ itself is invertible and
\eqref{eq:PLMEgen} is simplified to
\begin{equation}\label{eq:LMEnonsing}
	\lambda_J(x,y) = A_J^{-T}\nabla_y f(x,y).
\end{equation}
To show the differences between PLMEs and LMEs,
we give a toy example to compare their computational efficiency.

\begin{example}\label{ex:LMEdcp}
In \eqref{eq:genblv}, consider $n=1,p=2$ and
\[
	A = \left[\begin{array}{rr}
		1 & 0 \\1 & 0\\ 1 & 2\\ -1 & -1\end{array}\right],
	\quad b(x) = \bbm 0\\ 2-x\\ 2x^2\\-3\ebm.
\]
Since $m=4$ and $\rank\, A=2$, we have
\[
	\mc{P} = \{ \{1,3\},\, \{1,4\},\, \{2,3\},\, \{2,4\},\, \{3,4\}\}.
\]
The PLME in \eqref{eq:LMEnonsing} has degree $\deg(\nabla_y f)-1$
for each $J\in \mc{P}$. For instance,
\[
	\lambda_{\{1,3\}} = A_{\{1,3\}}^{-T}\nabla_y f
	= \bbm 1 & 0\\ 1 & 2\ebm^{-T} \bbm \nabla_{y_1} f\\ \nabla_{y_2} f\ebm
	=\bbm \nabla_{y_1} f\\ -\frac{1}{2}\nabla_{y_1}f+\frac{1}{2}\nabla_{y_2} f\ebm.
\]
Consider the LME \eqref{eq:lme} implied by $L(x,y)$ and $d(x,y)$ in \eqref{eq:Lnd}.
It is easy to compute the degree of $d(x,y)$ equals $10$ and
the degree of $L(x,y)\hat{f}(x,y)$ equals $\deg(\nabla_y f)+9$.
Therefore, PLMEs are much simpler than the rational LME for this problem.
\end{example}

\subsection{Decomposition of the KKT system}

With PLMEs, the KKT subset $\mc{K}_J$ can be explicitly expressed as
a semialgebraic set for every $J\in\mc{P}$.

\begin{prop}\label{prop:PLME}
For every $J\in \mc{P}$, let $\lambda_J(x,y)$ be as in \eqref{eq:PLMEgen}.
Then
\begin{equation}\label{eq:KKTsubset}
	\mc{K}_J = \left\{
	(x,y)\in \mc{U}\left|\begin{array}{l}
	A_J^T\lambda_J(x,y) = \nabla_y f(x,y),\\
	0\le [a_j^Ty-b_j(x)]\perp \lambda_j(x,y)\ge 0\,\,(j\in J)
	\end{array}\right.\right\}.
\end{equation}
In particular, if $\rank\, A = p$, then
\[
	\mc{K}_J = \left\{(x,y)\in \mc{U}\,\vert\,
	0\le [A_J^Ty-b_J(x)] \perp A_J^{-T}\nabla_y f(x,y)\ge 0\right\}.
\]
\end{prop}
\begin{proof}
For $J\in \mc{P}$, let $\hat{\mc{K}}_J$ denote the right-hand-side of \eqref{eq:KKTsubset}.
If $(u,v)\in \mc{K}_J$, then by \eqref{eq:defKJ},
there exists a vector of Lagrange multipliers $\lambda_J$ such that
$A_J^T\lambda_J = \nabla_y f(u,v)$.
Since $A_J^T$ has full column rank, we must have
\[
	\lambda_J = (A_JA_J^T)^{-1}A_J\nabla_y f(u,v) = \lambda_J(u,v).
\]
This implies $\mc{K}_J\subseteq \hat{\mc{K}}_J$.
Conversely, if $(u,v)\in\mc{\hat{K}}_J$,
then we also have $(u,v)\in \mc{K}_J$ with $\lambda_J = \lambda_J(u,v)$,
thus $\hat{\mc{K}}_J\subseteq\mc{K}_J$.
The rest conclusion is implied by \eqref{eq:KKTsubset}
\end{proof}

The KKT set $\mc{K}$ can be expressed through KKT subsets $\mc{K}_J$ as follows.

\begin{theorem}\label{thm:kktdcp}
Let $\mc{K}$ be the KKT set as in \eqref{eq:kktset}. Then
\begin{equation}\label{eq:KKTdcp3}
	\mc{K} \quad  = \bigcup_{  J\in \mc{P}} \mc{K}_J .
\end{equation}
\end{theorem}
\begin{proof}
The $\mc{K}\supseteq \bigcup_{J\in\mc{P}}\mc{K}_J$ since each
$\mc{K}_J\subseteq \mc{K}$.
It suffices to show the other direction.
Fix an arbitrary pair $(u, v) \in \mc{K}$ and denote
\[
	J_1 \,=\, \{ j \in [m]\,\vert \, a_j^Tv - b_j(u) = 0 \}.
\]
Let $r=\rank\ A_{J_1}$.
By Carath\'{e}odory's Theorem, there exist a subset
$J_2 \subseteq J_1$ and a Lagrange multiplier vector $\lambda = (\lambda_j)$
such that $\rank\, A_{J_2} = |J_2| = r$ and
\[
	\begin{array}{rl}
	\nabla_y f(u,v) - A_{J_2}^T\lambda_{J_2} = 0,\\
	0 \le [a_j^Tv-b_j(u)] \perp \lambda_j \ge 0\,  & (j \in J_2).
	\end{array}
\]
Since $r\le \rank\, A = t$, we can extend $J_2$ to some $J\in \mc{P}$
so that $(u,v)\in \mc{K}_J$.
As this argument works for arbitrary $(u,v)\in \mc{K}$, the conclusion holds.
\end{proof}

\begin{example}\label{ex:dcpvisual}
Consider the BPOP
\begin{equation*}\label{eq:locmin}
	\left\{\begin{array}{cl}
	\min\limits_{x,y\in \re} & (x-1.5)^2+y^2\\
	\st & x\ge 0,\,  2y+1\ge 0, \quad   y\in S(x),
	\end{array}\right.
\end{equation*}
where $S(x)$ is the set of optimizer(s) for
\[
	\left\{\begin{array}{cl}
	\min\limits_{y\in\re} & (y-x)^2\\
	\st & (y+1,\, 1-y,\, 4-2x-y,\, 3x-y-1)\ge 0.
	\end{array}\right.
\]
The $\mc{P} = \{\{1\},\{2\}, \{3\},\{4\}\}$ and
\[
	\mc{U} = \left\{(x,y)\in\re^2\,\Big\vert\,
	x\ge 0,\, -\frac{1}{2}\le y\le 1,\,
	4-2x-y\ge 0,\, 3x-y-1\ge 0 \right\}.
\]
By Proposition~\ref{prop:PLME}, the KKT subset $\mc{K}_{\{j\}}$
has an explicit expression for each $j\in[4]$,
which is presented in \Cref{tab:exdcpvisual}.
\begin{table}[htb]
\caption{PLMEs and KKT subsets for Example~\ref{ex:dcpvisual}}
\label{tab:exdcpvisual}
\centering
\begin{tabular}{|l|l|l|}
	\hline
	$J$ & $\lambda_J(x,y)$ & $\mc{K}_J$\\ [2pt] \hline
	$\{1\}$ & $ 2(y-x)$ & $\{(x,y)\in \mc{U}\,\vert\, 0\le (y+1)\perp2(y-x)\ge 0\}$\\  \hline
	$\{2\}$ & $2(x-y)$ & $\{(x,y)\in \mc{U}\,\vert\, 0\le (1-y)\perp 2(x-y)\ge 0\}$\\  \hline
	$\{3\}$ & $2(x-y)$ & $\{(x,y)\in \mc{U}\,\vert\, 0\le (4-3x-y)\perp 2(x-y)\ge 0\}$ \\  \hline
	$\{4\}$ & $2(x-y$) & $\{(x,y)\in \mc{U}\,\vert\, 0\le (3x-y-1)\perp 2(x-y)\ge 0\}$ \\  \hline
\end{tabular}
\end{table}
Since the lower-level problem is convex,
$\mc{K} = \mc{F}$ can be analytically solved as
\[
	\mc{K} \,=\, \Big \{\,  (x,y)\in\re^2\,\Big\vert\,
			\frac{1}{6}\le x\le \frac{9}{4},\,
			y = \min\{3x-1, x, 1, 4-2x\}  \,\Big \}.
\]
We plot these KKT sets in \Cref{fig:KKTdcp:a},
which shows the decomposition as in \eqref{eq:KKTdcp3}.
\begin{figure}[htbp]
\centering
\definecolor{ffqqqq}{rgb}{1,0,0}
\definecolor{yqqqqq}{rgb}{0.5,0,0}
\begin{tabular}{ccc}
	\begin{tikzpicture}[line cap=round,line join=round,>=triangle 45,x=1cm,y=1cm]
		\begin{axis}[
			x=1.1cm,y=1.1cm,
			axis lines=middle,
			xmin=-0.2,
			xmax=2.7,
			ymin=-1.2,
			ymax=1.5,
			xlabel = \(x\),
			ylabel = {\(y\)}]
			\draw [line width=1pt,color=ffqqqq] (1,1)-- (0.5,0.5);
			\draw [line width=1pt,color=ffqqqq] (1.5,1)-- (2.25,-1/2);
			\draw [line width=1pt,color=ffqqqq] (1,1)-- (1.5,1);
			\draw [line width=1pt,color=ffqqqq] (0.5,0.5)-- (1/6,-1/2);
		\end{axis}
	\end{tikzpicture} &
	\begin{tikzpicture}[line cap=round,line join=round,>=triangle 45,x=1cm,y=1cm]
		\begin{axis}[
			x=1.1cm,y=1.1cm,
			axis lines=middle,
			xmin=-0.2,
			xmax=2.7,
			ymin=-1.2,
			ymax=1.5,
			xlabel = \(x\),
			ylabel = {\(y\)}]
			\draw [line width=1pt,color=ffqqqq] (1,1)-- (0.5,0.5);
		\end{axis}
	\end{tikzpicture} &
	\begin{tikzpicture}[line cap=round,line join=round,>=triangle 45,x=1cm,y=1cm]
		\begin{axis}[
			x=1.1cm,y=1.1cm,
			axis lines=middle,
			xmin=-0.2,
			xmax=2.7,
			ymin=-1.2,
			ymax=1.5,
			xlabel = \(x\),
			ylabel = {\(y\)}]
			\draw [line width=1pt,color=ffqqqq] (1,1)-- (0.5,0.5);
			\draw [line width=1pt,color=ffqqqq] (1,1)-- (1.5,1);
		\end{axis}
	\end{tikzpicture} \\
	(a) $\mc{K}= \bigcup\limits_{j=1}^4 \mc{K}_{ \{j\}  }$
	& (b) $\mc{K}_{\{1\}}$ & (c) $\mc{K}_{\{2\}}$  \\[6pt]
	&
	\begin{tikzpicture}[line cap=round,line join=round,>=triangle 45,x=1cm,y=1cm]
		\begin{axis}[
			x=1.1cm,y=1.1cm,
			axis lines=middle,
			xmin=-0.2,
			xmax=2.7,
			ymin=-1.2,
			ymax=1.5,
			xlabel = \(x\),
			ylabel = {\(y\)}]
			\draw [line width=1pt,color=ffqqqq] (1,1)-- (0.5,0.5);
			\draw [line width=1pt,color=ffqqqq] (1.5,1)-- (2.25,-1/2);
		\end{axis}
	\end{tikzpicture} &
	\begin{tikzpicture}[line cap=round,line join=round,>=triangle 45,x=1cm,y=1cm]
		\begin{axis}[
			x=1.1cm,y=1.1cm,
			axis lines=middle,
			xmin=-0.2,
			xmax=2.7,
			ymin=-1.2,
			ymax=1.5,
			xlabel = \(x\),
			ylabel = {\(y\)}]
			\draw [line width=1pt,color=ffqqqq] (1,1)-- (0.5,0.5);
			\draw [line width=1pt,color=ffqqqq] (0.5,0.5)-- (1/6,-1/2);
		\end{axis}
	\end{tikzpicture} \\
	& (d) $\mc{K}_{\{3\}}$  & (e) $\mc{K}_{\{4\}}$  \\[6pt]
\end{tabular}
\caption{The PLME decomposition for the KKT set of Example~\ref{ex:dcpvisual}}
\label{fig:KKTdcp:a}
\end{figure}
\end{example}

\section{The Disjunctive Decomposition}\label{sc:disjunctive}

This section gives a disjunctive program relaxation of \eqref{eq:genblv}
based on \eqref{eq:KKTdcp3}.
For every $(x,y)\in\mc{K}$, we use
\begin{align}\label{set:Paxy}
	\mc{P}^{x,y} & \,\coloneqq\, \{J\in\mc{P}\,\vert\, (x,y)\in\mc{K}_J\},\\
	\label{eq:actset}
	\mc{I}(x,y) & \,\coloneqq\, \{ j\in[m] \,\vert\,  a_j^Ty = b_j(x) \}.
\end{align}
to denote the index set of branch KKT subsets where the point $(x,y)$ lies,
and the active index set of the lower-level inequalities, respectively.

\subsection{Disjunctive KKT relaxation}

Consider the relaxation of \eqref{eq:genblv}
\[
	(P): \, \left\{\begin{array}{cl}
			\min & F(x,y)\\
			\st & (x,y)\in \mc{K} = \bigcup\limits_{J\in\mc{P}} \mc{K}_J.
			\end{array}\right.
\]
For $J\in \mc{P}$, the $J$-th {\it branch problem} of $(P)$ is
\[
	(P_J):\, \left\{\begin{array}{rl}
		\vartheta_J\coloneqq \min & F(x,y)\\
			\st & (x,y)\in \mc{K}_J.
			\end{array}\right.
\]
Let $\vartheta_J$ denote the optimal value and
$(u_J,v_J)$ the optimizer of $(P_J)$, respectively.
It is clear that the smallest $\vartheta_J$ for $J\in \mc{P}$
is the optimal value of $(P)$ and the corresponding $(u_J,v_J)$
is an optimizer of $(P)$.
There are two major benefits for solving
$(P)$ through solving its branch problems $(P_J)$.
The first benefit is that $\mc{K}_J$ usually has much simpler
Lagrange multiplier expressions than the set $\mc{K}$. The second benefit is that
$(u_J,v_J)$ may be a local optimizer of the original BPOP, under some assumptions.
The local optimality can be verified by a convenient sufficient condition
using $\mc{I}(u_J,v_J)$.
To show the second benefit, we present a useful lemma.

\begin{lem}\label{lem:gloopt}
Suppose $(u,v)$ is a feasible point of \eqref{eq:genblv}.
	\begin{enumerate}
	
	\item[(i)]
	If $F(u,v)\le \vartheta_J$ for every $J\in\mc{P}^{u,v}$,
	then $(u,v)$ is a local optimizer of \eqref{eq:genblv}.
	
	\item[(ii)]
	If $F(u,v)\le \vartheta_J$ for every $J\in\mc{P}$,
	then $(u,v)$ is a global optimizer of \eqref{eq:genblv}.
	\end{enumerate}
\end{lem}
\begin{proof}
(i) Since $\mc{F}\subseteq \mc{K}$, $(u,v)$ is a local minimizer
of \eqref{eq:genblv} if there exists a scalar $\rho>0$ such that
$F(u,v)\le F(x,y)$ for all $(x,y)\in \mc{K}\cap B_{\rho}(u,v)$.
We first show that there exists $\rho>0$ that is sufficiently small
such that
\begin{equation}\label{Eqn3.3}
	\mc{K}\cap B_{\rho}(u,v) =
	\bigcup\limits_{J\in\mc{P}} \mc{K}_J\cap B_{\rho}(u,v)
	=\bigcup\limits_{J\in\mc{P}^{u,v}} \mc{K}_J\cap B_{\rho}(u,v).
\end{equation}
The first equality and the $\supseteq$ in the second equation obviously hold
for any $\rho>0$ by virtue of \cref{thm:kktdcp}.
It suffices to show the inclusion $\subseteq$ in the second equation.
For any $(u^{\rho}, v^{\rho})\in\mc{K}\cap B_{\rho}(u,v)$,
there exists $J\in\mc{P}$ such that $(u^{\rho}, v^{\rho})\in \mc{K}_J\cap B_{\rho}(u,v)$.
As $\rho\to 0$, $(u^{\rho}, v^{\rho})\to (u,v)$.
By the continuity of functions $\nabla_y f$, $b_j$ and $\lambda_j$,
we must have $(u,v)\in \mc{K}_J$ and $J\in \mc{P}^{u,v}$,
thus the second equality in \eqref{Eqn3.3} holds.
Since $\vartheta_J$ is the optimal value of $(P_J)$,
\eqref{Eqn3.3} implies that $F(u,v)\le F(x,y)$ for all
$(x,y)\in \mc{K}\cap B_{\rho}(u,v)$.
Therefore, $(u,v)$ is a local minimizer of \eqref{eq:genblv}.

(ii) By Theorem~\ref{thm:kktdcp}, if $F(x,v)\le \vartheta_J$ for all $J$,
then $F(u,v)\le F(x,y)$ for every $(x,y)\in \mc{K}$.
Note that $(u,v)\in \mc{F}\subseteq \mc{K}$.
So it is a global optimizer of \eqref{eq:genblv}.
\end{proof}

The formal results on the verification of
local optimality at $(u_J,v_J)$ for the original BPOP are summarized as follows.
\begin{theorem}\label{thm:eqv:PJ}
	Let $J\in \mc{P}$. Suppose $(P_J)$ has an optimizer $(u_J,v_J)$
	that is feasible for \eqref{eq:genblv}.
	\begin{enumerate}
	\item[(i)]
	If $\vartheta_J\le \vartheta_{J'}$ for every $J'\in \mc{P}$
	such that $J' \subseteq \mc{I}(u_J,v_J)\cup J$,
	then $(u_J, v_J)$ is a local optimizer of \eqref{eq:genblv}.
	
	\item[(ii)]
	If $\vartheta_J\le \vartheta_{J'}$ for all $J'\in\mc{P}$,
	then $(u_J,v_J)$ is a global optimizer of \eqref{eq:genblv}.
	\end{enumerate}
\end{theorem}
\begin{proof}
The result (ii) is directly implied by Lemma~\ref{lem:gloopt},
so it suffices to prove (i).

For every KKT point $(x,y)\in \mc{K}$,
let $\lambda^{x,y}$ be a corresponding Lagrange multiplier.
We can write the KKT equation as
\[
	\nabla_y f(x,y) = \sum\limits_{j\in\mc{I}(u_J,v_J)\cup J} \lambda_j^{x,y} a_j
	+\sum\limits_{j\not\in \mc{I}(u_J,v_J)\cup J}\lambda_j^{x,y} a_j.
\]
Since $b_j(x)$ is a polynomial, for all $(x,y)\in B_{\rho}(u_J,v_J)$,
$a_j^Tv_J>b_j(u_J)$ implies $a_j^Ty>b_j(x)$ when $\rho>0$ is sufficiently small.
In other words, there exists $\rho>0$ such that
\begin{equation}\label{eq:thm:eqvPJ}
	\mc{I}(x,y)\subseteq \mc{I}(u_J,v_J)\quad\forall (x,y)\in B_\rho(u_J, v_J).
\end{equation}
By the complementary slackness condition,
if $(x,y)\in \mc{K}\cap B_{\rho}(u_J,v_J)$, then
\[
	\lambda^{x,y}_j=0\quad \forall j \not \in \mc{I}(u_J, v_J).
\]
Therefore, for every $(x,y)\in \mc{K}\cap B_{\rho}(u_J,v_J)$, we have
\[
	\nabla_y f(x,y) = \sum\limits_{j\in\mc{I}(u_J,v_J)\cup J} \lambda_j^{x,y} a_j.
\]
By \Cref{thm:kktdcp} and Carath\'{e}odory's Theorem
(see e.g., \cite[Lemma 1]{AndHaeSchuSilva}),
the above equation implies $(x,y)\in\mc{K}_{J'}$ for some $J'\in \mc{P}$
that satisfies $J'\subseteq \mc{I}(u_J,v_J)\cup J$, thus
\begin{equation}\label{eq3.5}
	\mc{K}\cap B_{\rho}(u_J,v_J)
	\,\,=\,\, \bigcup_{J'\in \mc{P},J'\subseteq\mc{I}(u_J, v_J)\cup J}
	\mc{K}_{J'}\cap B_\rho (u_J,v_J).
\end{equation}
For every $J'\in \mc{P}$, since $\vartheta_{J'}$ is the optimal value
of $(P_{J'})$, we have
\[
	\vartheta_J = F(u_J, v_J) \ge \vartheta_{J'}\quad
	\mbox{if $J'\subseteq \mc{I}(u_J, v_J)\cup J$}.
\]
Under the assumption that $\vartheta_{J'}\ge \vartheta_J$ for all such $J'$,
$\vartheta_J$ must be equal to the minimum value of $f$ on
$\mc{K}\cap B_{\rho}(u_J,v_J)$.
So $(u_J,v_J)$ is a local minimizer of \eqref{eq:genblv}.
\end{proof}
Note that the condition of \Cref{thm:eqv:PJ}~(i) is weaker than
\Cref{lem:gloopt}~(i) since
\[
	\mc{P}^{u_J,v_J}\supseteq \{J'\in \mc{P}\,\vert\,
	J'\subseteq \mc{I}(u_J,v_J)\cup J\}.
\]
For the special case that $\mc{I}(u_J, v_J)\subseteq J$,
we have $J = \mc{I}(u_J,v_J)\cup J$.
The following result is simply implied by Theorem~\ref{thm:eqv:PJ}.
\begin{cor}\label{Cor3.3}
	Given $J\in \mc{P}$, suppose $(P_J)$ has an optimizer $(u_J,v_J)$
	that is feasible for \eqref{eq:genblv}.
	Then $(u_J,v_J)$ is a local optimizer of \eqref{eq:genblv}
	if $\mc{I}(u_J,v_J)\subseteq J$.
\end{cor}

Next, we use an example to illustrate the usage of \Cref{thm:eqv:PJ}.

\begin{example}\label{ex:locmin}
Consider the BPOP given in \Cref{ex:dcpvisual} where
\[
	\begin{array}{ll}
	F(x,y) = (x-1.5)^2+y^2, & f(x,y) = (y-x)^2,\\
	h(x,y) = (x,\, 2y+1), &
	g(x,y) = (y+1,\, 1-y,\, 4-2x-y,\, 3x-1-y).
	\end{array}
\]
Since the lower-level objective function $f$ is convex in $y$,
we have $\mc{F} = \mc{K}$ and each $\mc{K}_{\{j\}}\subseteq\mc{F}$,
whose parametric expressions are explicitly given in \Cref{ex:dcpvisual}.
For each $j\in[4]$, we solve $(P_{\{j\}})$ accurately with
all results reported in \Cref{tab:locmin}.
\Cref{fig:locmin} shows the feasible set along with the contour of
the upper-level objective function for this bilevel problem.
One can easily observe that the bilevel problem has a unique global optimizer
$(1.9,0.2)$ and a unique non-global local optimizer $(0.75,0.75)$.
\begin{table}[htb]
\caption{Optimizers for decomposed BPOPs of \Cref{eq:locmin}}
\label{tab:locmin}
\centering
\begin{tabular}{|c|c|l|c|c|}
	\hline
	$J$ & $\vartheta_J$ & $(u_J, v_J) $ & $\mc{I}(u_J,v_J)$ & $\mc{P}^{u_J,v_J}$\\ \hline
	$\{1\}$ & $1.125$ & $(0.75, 0.75)$ & $\emptyset$ & $\{1,2,3,4\}$\\ \hline
	$\{2\}$ & $1$ &  $(1.5,1)$ & $\{2, 3\}$ & \{2,3\}\\\hline
	$\{3\}$ & $0.2$ &    $(1.9, 0.2)$ & $\{3\}$ & \{3\} \\ \hline
	$\{4\}$ & $1.125$ & $(0.75,0.75)$ &  $\emptyset$ & $\{1,2,3,4\}$\\
	\hline
\end{tabular}
\end{table}

\begin{figure}[htb]
\centering
\definecolor{ffqqqq}{rgb}{1,0,0}
\definecolor{yqqqqq}{rgb}{0.5019607843137255,0,0}
\begin{tikzpicture}
	\begin{axis}[
		x=2.5cm, y=2.5cm,
		axis lines=middle,
		xmin= 0,
		xmax= 2.5,
		ymin=-0.7,
		ymax=1.2,
		xlabel = \(x\),
		ylabel = {\(y\)},
		legend pos=outer north east]
		\draw [line width=1pt,color=ffqqqq] (1,1)-- (0.5,0.5);
		\draw [line width=1pt,color=ffqqqq] (1,1)-- (1.5,1);
		\draw [line width=1pt,color=ffqqqq] (1.5, 1)-- (2.25,-0.5);
		\addplot[line width=1pt,color=ffqqqq] coordinates {(0.5,0.5) (1/6,-0.5)};
		\addlegendentry{feasible set $\mc{F}$};
		\addplot[domain = 0:2.5, samples=100,smooth, dashed, forget plot] {sqrt(1.125-(x-1.5)^2)} node[pos=0.7] {$1.125$};
		\addplot[domain = 0.7:2.3, samples=101,smooth, dashed, forget plot] {sqrt(0.5-(x-1.5)^2)} node[pos=0.3] { $0.6$ };
		\addplot[domain = 1:2, samples=100,smooth, dashed, forget plot] {sqrt(0.2-(x-1.5)^2)} node[pos = 0.4] { $0.2$ };
		\addplot[domain = 1:2, samples=100,smooth, dashed, forget plot] {-sqrt(0.2-(x-1.5)^2)};
		\addplot[domain = 0.7:2.3, samples=101,smooth, dashed, forget plot] {-sqrt(0.6-(x-1.5)^2)} ;
		\addplot[domain = 0:2.5, samples=100,smooth, dashed] {-sqrt(1.125-(x-1.5)^2)};
		\addlegendentry{contour of $F$};
		\addplot [only marks, mark = square, thick] coordinates{(1.9,0.2)};
		\addlegendentry{global min};
		\addplot [only marks, mark = star, thick] coordinates{(0.75,0.75)};
		\addlegendentry{local min};
	\end{axis}
\end{tikzpicture}
\caption{For Example~\ref{ex:locmin},
	the solid line is the feasible set $\mc{F}$,
	the dotted curves are contour level curves of $F$,
	the ``$\Box$'' denotes the global minimizer
	and the ``$\star$'' denotes the local minimizer.
}
\label{fig:locmin}
\end{figure}
Subsequently, we discuss the optimality verification for each computed $(u_J,v_J)$.

For $J = \{3\}$, we have $(u_{J},v_{J}) = (1.9,0.2),\,\vartheta_J=0.2$
from the table.
From the graph, $ (1.9,0.2)\in \mathcal{F}$.
By \Cref{thm:eqv:PJ}, since $ \vartheta_J=0.2 \leq  \vartheta_{J'}$
for all $J'\in \mc{P}$,
we detect $(x,y) = (1.9,0.2)$ as the global optimizer.

For $J=\{1\}, \{4\}$, we have
$(u_{J},v_{ J} )=(0.75,0.75),\, \mc{I}(u_J,v_J)=\emptyset$
from the table.
From the graph, $ (0.75,0.75)\in \mathcal{F}$.
By \Cref{Cor3.3}, we verify that $(x,y) =(0.75,0.75) $ is a local
minimizer of this bilevel problem.

For $J = \{2\}$, we have $(u_{J},v_{ J})=(1.5,1),\,
\mc{I}(u_J,v_J)=\{2,3\}$ from the table.
From the graph, $(1.5, 1)\in \mathcal{F}$.
Since $\mc{I}(u_J,v_J) \cup J=\{2,3\}$ and
$ \vartheta_{\{3\}}=0.2<  \vartheta_{\{2\}}$,
we cannot decide the local optimality of $(x,y) = (1.5,1)$
by using Theorem~\ref{thm:eqv:PJ}.
But this is correct since $(1.5,1)$ is not a local optimum.
\end{example}

\subsection{Local optimality of BPOPs}\label{ssc:lc}

We would like to remark that Theorem~\ref{thm:eqv:PJ}
does not give a full characterization for local optimizers.
In computational practice, local optimality can be further verified
by solving branch problems with an additional small ball constraint.
Let $(\hat{u},\hat{v})\in\mc{F}$ denote the selected point.
The verification process is as follows.	
First, determine the active set $\mc{I}(\hat{u}, \hat{v})$
and a label set $\hat{J}\in\mc{P}$
such that $(\hat{u},\hat{v})\in \mc{K}_{\hat{J}}$.
Then select a small scalar $\rho>0$ and solve
\begin{equation}\label{eq:locminver}
	(Q_J):\, \left\{\begin{array}{rl}
			\varphi_J\coloneqq \min & F(x,y) - F(\hat u,\hat v)\\
			\st  & (x,y)\in \mc{K}_{J}\cap B_\rho(\hat u, \hat v),
		\end{array}\right.
\end{equation}
for every $J\in \mc{P}$ with $J\subseteq \mc{I}(\hat{u}, \hat{v})\cup \hat{J}$.
Let $\varphi_J$ denote the optimal value of $(Q_J)$.
If each $\varphi_J\ge 0$, then $(\hat{u},\hat{v})$ is a local minimizer
of \eqref{eq:genblv}.
We summarize this result in the following theorem.
\begin{theorem}\label{thm:loc}
Let $(\hat{u},\hat{v})$ be a feasible point of \eqref{eq:genblv}
and let $\rho>0$ be sufficiently small.
Suppose $(\hat{u},\hat{v})\in \mc{K}_{\hat{J}}$ for some $\hat{J}\in \mc{P}$.
If $\varphi_J\ge 0$ for all $J\in \mc{P}$ such that
$J\subseteq \mc{I}(\hat{u},\hat{v})\cup \hat{J}$, then
$(\hat{u},\hat{v})$ is a local minimizer of \eqref{eq:genblv} .
\end{theorem}
\begin{proof}
Similar to the proof of \eqref{Eqn3.3} and \eqref{eq3.5},
we can show that
for $\rho>0$ that is sufficiently small,
it holds
\[
	\mc{F}\cap B_{\rho}(\hat{u}, \hat{v})\subseteq
	\bigcup_{J\in\mc{P},J\subseteq \mc{I}(\hat{u}, \hat{v})\cup \hat{J}}
	\mc{K}_J\cap B_{\rho}(\hat{u}, \hat{v}).
\]
Since $\varphi_J$ is the optimal value of $(Q_J)$,
for every $(x,y)\in \mc{K}_J\cup B_{\rho}(\hat{u}, \hat{v})$,
we have
\[
	F(x,y) \, =\,F(\hat{u}, \hat{v}) + \big( F(x,y) -F(\hat{u}, \hat{v}) \big)
			\ge \,F(\hat{u}, \hat{v})+ \varphi_J.
\]
If $\varphi_J\ge 0$ for all $J\in \mc{P}$ such that
$J\subseteq \mc{I}(\hat{u},\hat{v})\cup \hat{J}$, then
\[
	F(x,y)\ge F(\hat{u},\hat{v}),\quad
	\forall (x,y)\in \mc{F}\cap B_{\rho}(\hat{u}, \hat{v}).
\]
It follows that $(\hat{u}, \hat{v})$ is a local minimizer of \eqref{eq:genblv}.
\end{proof}

We remark that \eqref{eq:locminver} is a polynomial optimization problem,
which can be solved globally by Moment-SOS relaxations.
To verify that $(\hat{u}, \hat{v})$ is a local minimizer of \eqref{eq:genblv},
it is enough to choose a priori value for $\rho>0$ (say, $0.1$)
such that the optimal value $\varphi_J$ of \eqref{eq:locminver} is nonnegative
for all $J\in \mc{P}$ with 
$J\subseteq \mc{I}(\hat{u}, \hat{v}) \cup \hat{J}$ 
(given $(\hat{u}, \hat{v})\in \mc{K}_{\hat{J}}$).
It is pretty hard to get a simple lower bound for such $\rho$.
In fact, there does not exist a uniform lower bound for $\rho$
in polynomial optimization \cite{nie2015local}.
However, in computational practice, the local optimality can be ensured
if we get one value for $\rho$ that satisfies all conditions of \cref{thm:loc}.
By this observation, we give \cref{def:alg:locmin}
to verify local optimality of $(\hat{u}, \hat{v})$.
It is interesting to note that the task of
verifying local optimality is NP-hard \cite{MurKab87}.

\section{The PLME-Disjunctive Algorithm}\label{sec:alglocal}

In this section, we give \Cref{def:alg:bpop}
to solve the BPOP \eqref{eq:genblv}.
For $J\in\mc{P}$, define the $J$th branch problem of \eqref{eq:genblv} as
\begin{equation}\label{eq:PJstar}
	(P_J^*):\,\left\{\begin{array}{rl}
		\vartheta_J^*\coloneqq \min & F(x,y)\\
		\st & (x,y)\in\mc{K}_J\cap \mc{F}.
		\end{array}\right.
\end{equation}
We assume each $(P_J^*)$ is either infeasible or has a global optimizer.
Under some general assumptions, the branch problem $(P_J^*)$
can be solved globally by \cref{def:alg:j},
which is called in an inner loop of \cref{def:alg:bpop}.

\begin{algorithm}
\caption{The main algorithm to solve \eqref{eq:genblv}.}
\label{def:alg:bpop}
\begin{algorithmic}[1]
\STATE{Let $\mc{J}_0, \mc{J}_1, \mc{J}_{2}$ be empty sets.}
\FOR{$J\in\mc{P}$}
\STATE{Solve $(P_J^*)$ with {\bf\Cref{def:alg:j}}.}
\IF{$(P_J^*)$ is solvable}
\STATE{Get the optimal value $\vartheta_J^*$
	and an optimizer $(u_J^*, v_J^*)$ of $(P_J^*)$.}

{Update $\mc{J}_0 \coloneqq \mc{J}_0\cup J$.}
\ENDIF
\ENDFOR
\IF{$\mc{J}_0 = \emptyset$}
\RETURN {The BPOP \eqref{eq:genblv} is infeasible.}
\ELSE
\STATE{Compute $F_{min}\coloneqq \min\{\vartheta_J^*\,\vert\,J\in\mc{J}_0\}$
	and $\mc{J}_1\coloneqq \{J\in \mc{J}_0\,\vert\, F_{min} = \vartheta_J^*\}$.}
\FOR{$J\in \mc{J}_0\setminus \mc{J}_1$}
\IF{$\vartheta_J^*\le \vartheta_{J'}^*$ for all $J'\in \mc{J}_0$
	such that $J'\subseteq \mc{I}(u_J^*,v_J^*)\cup J$}
\STATE Update $\mc{J}_{2} \coloneqq \mc{J}_2\cup J.$
\ENDIF
\ENDFOR
\RETURN $F_{min}$, $\{\vartheta_J^*, (u_J^*,v_J^*)\,\vert\, J\in\mc{J}_0\}$,
and $\mc{J}_0$, $\mc{J}_1,\mc{J}_2$.
\ENDIF
\end{algorithmic}
\end{algorithm}
In \Cref{def:alg:bpop}, the symbol $\mc{J}_0$ is used to store the index of
\eqref{eq:PJstar} that is feasible,
$\mc{J}_1$ is used to store the index of \eqref{eq:PJstar} whose optimizer(s)
are global optimizers of \eqref{eq:genblv}
and $\mc{J}_2$ is used to store the index of \eqref{eq:PJstar} whose optimizer(s)
are verified as non-global local optimizers of \eqref{eq:genblv}.

\begin{theorem}
Suppose all branch problems are solved globally by \Cref{def:alg:bpop}.
Then $F_{min}$ is the optimal value of \eqref{eq:genblv}.
For $J\in\mc{J}_0$, $(u_J^*, v_J^*)$ is a global optimizer of \eqref{eq:genblv}
if and only if $J\in\mc{J}_1$,
and is a non-global local optimizer of \eqref{eq:genblv} if $J\in \mc{J}_2$.
\end{theorem}
\begin{proof}
These conclusions are implied by \Cref{thm:kktdcp} and \Cref{thm:eqv:PJ}.
\end{proof}
We remark that $(u_J^*, v_J^*)$ may be a local optimizer of \eqref{eq:genblv}
even if $J\not\in \mc{J}_1\cup \mc{J}_2$.
The local optimality can be verified by using the method described in \Cref{ssc:lc}.
We summarize this as \cref{def:alg:locmin}, where the bisection trick
is used to get a priori value of $\rho$ such that the condition of \cref{thm:loc} is satisfied.
We remark that such $\rho$ always exists if $(u_J^*, v_J^*)$ is a local optimizer.
However, in the worst case, such value of $\rho$ can be very small,
which may cause numerical issues.
The problem of verifying local optimality is NP-hard \cite{MurKab87}.

\begin{algorithm}
\caption{The algorithm to verify local optimality of $(\hat{u},\hat{v})$.}
\label{def:alg:locmin}
\begin{algorithmic}[1]
	\STATE{Let $\rho = 0.1$. For the given candidate solution
		$(\hat{u},\hat{v})\in \mc{K}_{\hat{J}}\cap \mc{F}$, determine the index set
		 $\mc{J}' = \{J\in \mc{J}_0
		\,\vert\, J\subseteq \mc{I}(\hat{u},\hat{v})\cup \hat{J},\, \vartheta_{J}^*< \vartheta_{\hat{J}}^*\}$.}
	\WHILE{ $\mc{J}'\not= \emptyset$ }
		\FOR{$J\in \mc{J}'$}
				\STATE Solve \eqref{eq:locminver} to get an optimal value $\varphi_{J}$.\\
				\IF{$\varphi_{J}\ge 0$}
				\STATE Update $\mc{J}'\coloneqq \mc{J}'\setminus J$.
				\ELSE
				\STATE Update $\rho\coloneqq 0.5\rho$.
				\ENDIF
		\ENDFOR
		\IF{$\mc{J}'=\emptyset$}
			\RETURN{$(\hat{u}, \hat{v})$ is verified as a local optimizer with the radius $\rho$.}
		\ENDIF	
	\ENDWHILE
\end{algorithmic}
\end{algorithm}
In \cref{def:alg:locmin}, the index set $\mc{J}_0$ and optimal values
$\vartheta_{\hat{J}}^*, \vartheta_J^*$ are computed from \cref{def:alg:bpop},
and $(\hat{u}, \hat{v}) = (u_{\hat{J}}^*,v_{\hat{J}}^*)$
is a candidate solution obtained from some
branch problem $(P_{\hat{J}}^*)$ with $\hat{J}\not\in \mc{J}_1\cup \mc{J}_2$.
For convenience, we set the initial value of $\rho$ as $0.1$,
but it can be changed to any other small scalar.
By arguments in \cref{thm:loc}, when $\rho$ is sufficiently small,
every $(u,v)\in \mc{F}\cap B_{\rho}(\hat{u}, \hat{v})$
must belong to some $\mc{K}_{J}$ over
   $\{J\in \mc{J}_0\vert J\subseteq \mc{I}(\hat{u}, \hat{v}) \cup \hat{J}\}$.
Fix such an index set $J$.
Since $\vartheta_{J}^*$ denotes the optimal value of $(P_J^*)$,
if $\vartheta_{J}^*\ge \vartheta_{\hat{J}}^*$, then the optimal value of \eqref{eq:locminver}
satisfies $\varphi_{J}\ge \vartheta_{J}^*-\vartheta_{\hat{J}}^*\ge 0$;
otherwise, we need to compute $\vartheta_J$ by solving \eqref{eq:locminver} with the given $\rho$.
This further verification process is made for all $J\subseteq \mc{I}(\hat{u}, \hat{v})\cup \hat{J}$
with $\vartheta_{J}^*< \vartheta_{\hat{J}}^*$.
If all such $\varphi_{J}\ge 0$, then such $\rho$
is a valid radius to certify the local
optimality of  $(\hat{u}, \hat{v})$.
If some $\vartheta_J<0$, then we need to update $\rho$ with the bisection trick and
repeat the process until the termination criterion is met.

\begin{theorem}
If \cref{def:alg:locmin} terminates within finitely many loops,
then $\rho$ in the termination loop satisfies all conditions in \cref{thm:loc},
and thus $(\hat{u}, \hat{v})$ is a local optimizer of \eqref{eq:genblv}.
\end{theorem}
\begin{proof}
	This conclusion is implied by \cref{thm:loc}.
\end{proof}

\subsection{Solving $(P_J^*)$}\label{ssc:alg:blvsub}

We discuss how to solve the branch problem $(P_J^*)$.
In general, $(P_J^*)$ is easier to solve than
\eqref{eq:genblv}. However, in the worst case,
$(P_J^*)$ may be as hard as \eqref{eq:genblv}.
For instance, this is the case if $J=[m]$.
We make the following assumption for $(P_J^*)$.

\begin{ass}\label{ass:FE}
For given $J\in\mc{P}$ and $(\hat{x}, \hat{y})\in\mc{K}_J$, $\hat{z}\in S(\hat{x})$,
there exists a polynomial vector $q: \re^{n}\times \re^{p}\to \re^p$
such that
\begin{equation}\label{eq:FE:Kj}
	q(\hat{x}, \hat{y}) = \hat{z}\quad\mbox{and}\quad
	q(x,y)\in Z(x)\,\, \forall (x,y)\in \mc{K}_J,
\end{equation}
where $Z(x)$ is the lower-level feasible set of \eqref{eq:genblv}.
\end{ass}

Such a $q$ is called a {\it feasible extension} of $(\hat{x}, \hat{y})$
at $\hat{z}$ on $\mc{K}_J$.
Feasible extensions can be used to exclude infeasible points obtained from relaxations.
Under Assumption~\ref{ass:FE}, consider the inequality constraint
\begin{equation}\label{ineq:fe}
	f(x,q(x,y)) - f(x,y)\ge 0.
\end{equation}
If $\hat{y}\not\in S(\hat{x})$, then \cref{ineq:fe}
is violated at $(\hat{x}, \hat{y})$ because
\[
	f(\hat{x}, q(\hat{x}, \hat{y})) \,=\, f(\hat{x}, \hat{z})
	\,=\,  \inf\limits_{y\in Z(\hat{x})} f(\hat{x}, y) \,<\, f(\hat{x}, \hat{y}).
\]
However, every feasible point of \cref{eq:genblv} must satisfy \cref{ineq:fe}.
So, we can use feasible extensions and \cref{ineq:fe} to exclude infeasible points.
Feasible extensions can be obtained for many optimization problems.
They always exist if $\mc{K}_J$ is a finite set \cite[Theorem~4.2]{NieTangZgnep21}.
In particular, there are universal feasible extensions for single bound,
box and simplex constraints \cite{nie2020bilevel}.
We list some of them as follows.
For given $J\in \mc{P}$, fix $(\hat{x}, \hat{y})\in \mc{K}_J$ and $\hat{z}\in S(\hat{x})$.
Let $Z(x)$ denote the lower-level feasible set of \eqref{eq:genblv}.
\begin{itemize}

\item Suppose $Z(x) = \{y\in \re^p\,\vert\, y\ge l(x)\}$,
where $l = (l_1,\ldots, l_p)$ is a given polynomial vector.
Then $q(x,y) = l(x)+(\hat{z}-l(\hat{x}))$ satisfies Assumption~\ref{ass:FE}.

\item Suppose $Z(x) = \{y\in \re^p\,\vert\, l(x)\le y\le w(x)\}$,
where $l = (l_1,\ldots, l_p)$ and $w = (w_1,\ldots, w_p)$ are given polynomial tuples.
Then  $q = (q_1,\ldots, q_p)$ with
\[
	q_i(x,y) = \frac{w_i(\hat{x})-\hat{z}_i}{w_i(\hat{x})-l_i(\hat{x})}l_i(x)+
	\frac{\hat{z}_i-l_i(\hat{x})}{w_i(\hat{x})-l_i(\hat{x})}u_i(x)
\]
for each $i\in[p]$ satisfies Assumption~\ref{ass:FE}.

\item Suppose $Z(x) = \{y\in\re^p\,\vert\, y\ge l(x),\, w(x)-e^Ty\ge 0\}$,
where $w\in\re[x]$ and $l = (l_1,\ldots, l_p)$ is a given polynomial vector.
Then $q = (q_1,\ldots, q_p)$ with
\[
	q_i(x,y) = l_i(x)+\frac{\hat{z}_i-l_i(\hat{x})}{w(\hat{x})-e^Tl(\hat{x})}(u(x)-e^Tl(x))
\]
for each $i\in[p]$ satisfies Assumption~\ref{ass:FE}.
\end{itemize}
For general linear constraints, it is difficult to analytically check whether
or not Assumption~\ref{ass:FE} is satisfied.
However, if they exist, we can compute feasible extensions for general linear
and quadratic constraints by using the efficient method introduced in \Cref{sc:FE}.
It is an interesting question to find convenient conditions
for Assumption~\ref{ass:FE} to hold.
We plan to explore it as future work.

Summarizing the above, we give \cref{def:alg:j} to solve
the branch problem $(P_J^*)$ for each $J\in\mc{P}$.

\begin{algorithm}
\caption{The algorithm to solve $(P_J^*)$.}
\label{def:alg:j}
\begin{algorithmic}[1]
\STATE{Let $\mc{Q}_0 \coloneqq \emptyset$, $k\coloneqq 0$, $\eta_0 \coloneqq -1$.}
\WHILE{$\eta_k<0$}
\STATE{Solve the optimization problem
	\begin{equation}\label{eq:Kj:alg}
		\left\{\begin{array}{rl}
			\vartheta_k\coloneqq \min & F(x,y)\\
			\st & f(x,q(x,y))-f(x,y)\ge 0,\,\forall q\in \mc{Q}_k,\\
			& (x,y)\in\mc{K}_J.
		\end{array}\right.
\end{equation}}
\IF{\eqref{eq:Kj:alg} has an optimizer $(x^{(k)}, y^{(k)})$}
\STATE{Verify the lower-level feasibility of $(x^{(k)}, y^{(k)})$ by solving
	\begin{equation}\label{eq:feascheck}
		\left\{\begin{array}{rl}
			\eta_k \coloneqq \min & f(x^{(k)},y)-f(x^{(k)},y^{(k)})\\
			\st\, & Ay-b(x^{(k)})\ge 0.
		\end{array}\right.
\end{equation}}
\IF{$\eta_k\ge 0$}
\RETURN $\vartheta_J^* \coloneqq \vartheta_k$,
$u_J^*\coloneqq x^{(k)}$, $v_J^*\coloneqq y^{(k)}$.
\ELSE
\STATE{Solve for an optimizer $z^{(k)}$ of \cref{eq:feascheck}}.\\
Compute a feasible extension $q^{(k)}$ of $(x^{(k)}, y^{(k)})$
at $z^{(k)}$ on $\mc{K}_J$.\\
Update $\mc{Q}_{k+1}\coloneqq \mc{Q}_k\cup \{q^{(k)}\}$ and $k\coloneqq k+1$.
\ENDIF
\ELSE
\RETURN $(P_J^*)$ is infeasible.
\ENDIF
\ENDWHILE
\end{algorithmic}
\end{algorithm}
The convergence of \cref{def:alg:j} is shown in \Cref{ssc:convprop}.
Under Assumption~\ref{ass:FE}, if \cref{def:alg:j} terminates within finitely many loops,
then it either detects the infeasibility of $(P_J^*)$ or returns its optimal value
and optimizer.
This is the situation observed in all of our numerical experiments.
In fact, when $f$ is convex in $y$,
\cref{def:alg:j} always terminates at the initial loop $k=0$
(see Proposition~\ref{prop:conv:alg}).
However, it is possible that \cref{def:alg:j}
does not terminate within finitely many loops
in the worst case. In \cref{thm:asyconv:algj},
we show that the computed sequence of $(x^{(k)}, y^{(k)})$
will converge to a global optimizer of $(P_J^*)$ under some general assumptions.

\begin{figure}[!ht]
	\centering
	\begin{tikzpicture}[
		>=Stealth,
		box/.style={
			draw,
			rectangle,
			minimum width=1.2cm,
			minimum height=0.4cm, 
			align=center 
		},
		label/.style={font=\small, align=center}, 
		]

		\node[box] (P11) at (2.0, 5) {(1.1)};
		\node[box] (P) at (5, 5) {(P)};
		\node[box] (PJ) at (8.5, 5) {$(P_J)$};
		\node[box] (P44) at (8.5, 3) {(4.4)};
		\node[box] (PJs) at (2.0, 3) {$(P_J^*)$};

		\draw [->] (P11.east) -- (P.west) node[midway, above, label, yshift=2pt] {relaxation};

		\draw [->] (P.east) -- (PJ.west) node[midway, above, label, yshift=2pt] {decomposition};

		\draw [->] (PJ.south) -- (P44.north) node[midway, right=2pt, label] {add constraints with \\ feasible extensions};

		\draw [->] (P44.west) -- (PJs.east) node[midway, below, label, yshift=-2pt] {check lower-level feasibility by (4.5)};

		\draw [->] (PJs.north) -- (P11.south) node[midway, left=2pt, label] {check local/global \\ optimality};
		
	\end{tikzpicture}
	\caption{The diagram for relationships between branch problems}
	\label{fig:diagram}
\end{figure}
The diagram in \Cref{fig:diagram} shows the relationship
between involving problems in \Cref{def:alg:bpop,def:alg:j}.
The original bilevel problem \eqref{eq:genblv} is hard to solve directly.
We consider the relaxation problem $(P)$ obtained by replacing the optimality constraint
$y\in S(x)$ with the KKT set constraint $(x,y)\in \mc{K}$.
By PLMEs, $\mc{K}$ can be written as a union of subsets $\mc{K}_J$ such that
each $\mc{K}_J$ has a convenient expression.
This motivates us to solve $(P)$ by its branch problems $(P_J)$ with $J\in \mc{P}$.
When $f$ is convex, every KKT point is feasible for \eqref{eq:genblv},
thus $(P)$ is always a tight relaxation.
When $f$ is nonconvex, the optimizer of $(P_J)$ may not be feasible for the original
bilevel problem. For this case, we need to iteratively tighten $(P_J)$
with feasible extensions and check the lower-level feasibility of computed solutions,
until an optimizer of $(P_J^*)$ is obtained.
The local/global optimality of computed feasible solutions can be verified by \cref{thm:eqv:PJ}.

Each optimization problem in \Cref{def:alg:j} is polynomial optimization,
which can be solved globally by Moment-SOS relaxations.
Consider a general polynomial optimization problem
($z$ can be changed to other variables)
\begin{equation}\label{ppp}
	\left\{\begin{array}{rl}
		c^*\coloneqq \min\limits_{z\in\re^p} &  c_0(z)\\
		\st & c_i(z)\geq 0\,(i\in\mc{I}),
	\end{array}\right.
\end{equation}
where $c_0, c_i\,(i\in\mc{I})$ are polynomials and
$\mc{I}$ is a finite index set. Let
\[
	d_c\coloneqq \max\limits_{i\in\mc{I}}
	\{1,\lceil \deg(c_i)/2\rceil \},\quad
	d_0 \coloneqq \max\,\{d_c,\lceil \deg(c_0)/2\rceil\},
\]
where $\deg(\cdot)$ denotes the degree of a polynomial.
For an integer $l\ge d_0$, let $w = (w_{\alpha})\in\re^{\N_{2l}^p}$
be a real vector labeled with
\[
	\N_{2l}^p\coloneqq\{\alpha = (\alpha_1,\ldots, \alpha_p)\in\N^p\,\vert\,
	e^T\alpha\le 2l\}.
\]
We can define the bilinear operator
\[
	\langle c,w\rangle \coloneqq
	\sum\limits_{\alpha\in\N_{2l}^p} c_{\alpha}w_{\alpha}
	\quad \mbox{for every}\quad c(z) =
	\sum\limits_{\alpha\in\N_{2l}^p} c_{\alpha}z^{\alpha}.
\]
If there exists $\hat{z}\in\re^p$ such that $w_{\alpha}= \hat{z}^{\alpha}$
for every $\alpha\in \N_{2l}^p$, then $\langle c, w\rangle = c(\hat{z})$ for every
polynomial $c(z)$ with $\deg(c)\le 2l$.
Necessary conditions for this to hold imply the
$l$-th order moment relaxation of \eqref{ppp}:
\begin{equation}\label{dpr}
	\left\{\begin{array}{rl}
		c^{(l)}\coloneqq &  \langle c_0,w\rangle\\
		\st & L_{c_i}^{(l)}[w]\succeq 0\,(i\in\mc{I}),\\
		& M_{l}[w]\succeq 0,\\
		& w_0 = 1,\, w\in \re^{\N_{2l}^p},
	\end{array}\right.
\end{equation}
where $L_{c_i}^{(l)}[w]$ the $l$-th order localizing matrix of $c_i$ generated by $w$
and $M_l[w]$ is the $l$-th order moment matrix of $w$;
see \cite[Subsection~2.5]{MPO} for the notation.
The \eqref{dpr} is a semidefinite program.
Let $c^*,c^{(l)}$ denote the optimal values of \eqref{ppp} and \eqref{dpr}
respectively. Under archimedeanness,
$c^{(l)}\rightarrow c^*$ as $l\rightarrow\infty$ \cite{Las01}.
In addition, if some classical optimality conditions are satisfied,
then $c^{(l)} = c^*$ for all $l$ big enough,
as shown in \cite{nie2014optimality}.
Suppose $w^*$ is an optimizer of \eqref{dpr}.
The finite convergence can be conveniently checked
by the {\it flat truncation} condition:
there exists an integer $d_1\in [d_0,l]$ such that
\begin{equation}\label{eq:flat}
	r = \rank\, M_{d_1-d_c}[w^*] \,=\, \rank\, M_{d_1}[w^*].
\end{equation}
If \eqref{eq:flat} is satisfied, then $c^{(l)} = c^*$ and there exist $r$ optimizers
$v_1,\ldots, v_r$ of \eqref{ppp} and positive scalars $\mu_1,\ldots, \mu_r$ such that
\[
	w^*|_{2d_1} = \mu_1[v_1]_{2d_1}+\cdots +\mu_r[v_r]_{2d_1}.
\]
In the above, $w^*|_{2d_1}\coloneqq (w_{\alpha}^*)_{\alpha\in\N_{2d_1}^p}$ and
(for $z = v_1,\ldots, v_r$)
\[
	[z]_{2d_1} \coloneqq \bbm 1 & z_1 & \cdots & z_p & z_1^2 & z_1z_2 &
\cdots & z_p^{2d_1}\ebm^T.
\]
We refer to \cite[Lemma~5.3.1]{MPO} for the proof of this result.
These optimizers $v_1,\ldots, v_r$ can be computed with
Schur decompositions \cite{HL05},
the implementation of which is embedded in the software package
{\tt GloptiPoly~3} \cite{GloPol3}.
In our numerical experiments, polynomial optimization problems are all
solved globally under the flat truncation.

In \cref{def:alg:j}, the optimization \cref{eq:feascheck} is used to
verify the lower-level feasibility of $(x^{(k)}, y^{(k)})$.
This step can be skipped for the special case that
$f$ is convex in $y$; see the related result in Proposition~\ref{prop:conv:alg}.
The feasible extension $q^{(k)}$ satisfies
\begin{equation}\label{eq:fe:alg}
	q^{(k)}(x^{(k)}, y^{(k)}) = z^{(k)},\quad
	q^{(k)}(x,y)\in Z(x)\,\forall (x,y)\in \mc{K}_J.
\end{equation}
For each $y^{(k)}$ that does not optimize \eqref{eq:feascheck},
the inequality constraint
\begin{equation}\label{eq:ineq}
	f(x,q^{(k)}(x,y)) - f(x,y)\ge 0
\end{equation}
prevents the corresponding $(x^{(k)}, y^{(k)})$ to be feasible in future loops.
But the same constraint is satisfied at every feasible point of $(P_J^*)$,
so \eqref{eq:Kj:alg} is always a relaxation of $(P_J^*)$.
Note that \eqref{eq:ineq} may exclude several KKT points
that are not feasible for \eqref{eq:genblv}.
Consider the special case that \eqref{eq:ineq} is violated at every
$(x,y)\in\mc{K}_J$ at some loop $k$.
Then $(P_J^*)$ is infeasible and \eqref{eq:Kj:alg} is easy to solve at the same loop.
When $k$ is small, the computational time of \eqref{eq:Kj:alg} in each loop is similar.
In our numerical experiments, \Cref{def:alg:j} usually converges quickly
and terminates within a finite number of iterations.
However, for the worst case that \Cref{def:alg:j} does not terminate in finite loops,
\eqref{eq:Kj:alg} will become very challenging to solve with the increase of constraints.

The framework for updating relaxations with additional constraints is called
the {\it exchange method} or {\it cutting plane/surface method} in
semi-infinite program \cite{HettKor93}.
Many numerical algorithms for bilevel optimization are built in
combination with exchange methods such as
the bounding algorithm in \cite{mitsos2008global},
semidefinite algorithms in \cite{nie2017bilevel,nie2020bilevel} and
bound-and-cut algorithms in \cite{FisLju18,KleinerLa21,TaRaDe20}.
A common strategy is to approximate \eqref{eq:genblv} by
\[
	\left\{\begin{array}{cl}
	\min\limits_{(x,y)\in\mc{U}} & F(x,y)\\
	\st & v(x)-f(x,y)\ge 0,\quad\forall y\in Z,
	\end{array}\right.
\]
where $Z$ is a finite set and $v(x)$ is the value function of
the lower-level problem
\begin{equation}\label{eq:valfun}
	v(x) \, \coloneqq \, \inf_{z\in Z(x)} f(x,z).
\end{equation}
In \Cref{def:alg:j}, we apply the exchange method to \eqref{eq:Kj:alg},
whose feasible set is the KKT subset $\mc{K}_J$ at the initial loop $k=0$.
Since $\mc{K}$ has a smaller dimension than $\mc{U}$,
it is much easier for \eqref{eq:Kj:alg} to approximate $(P_J^*)$
by the exchange method.
We refer to \cite[Example~3.2]{nie2017bilevel} for the computational advantage
of applying the exchange method in combination of KKT set representations.
Under Assumption~\ref{ass:FE},
the problem~\eqref{eq:Kj:alg} is a relaxation of $(P_J^*)$.

\subsection{Convergence property}\label{ssc:convprop}

For the case of convex lower-level optimization,
\cref{def:alg:j} terminates at the initial loop.

\begin{prop}\label{prop:conv:alg}
	Assume $f(x,y)$ is convex in $y$ for every $(x,y)\in \mc{U}$.
	Then for each $J\in\mc{P}$,
	the optimization $(P_J)$ is equivalent to $(P_J^*)$ and \cref{def:alg:j}
	terminates at the initial loop $k=0$.
\end{prop}
\begin{proof}
By given conditions, every KKT point is equivalent to an optimizer,
i.e., $\mc{K} = \mc{F}$.
Then $\mc{K}_J = \mc{F}\cap \mc{K}_J$ for every $J\in\mc{P}$.
It implies that $(P_J)$ and $(P_J^*)$ are equivalent.
So \Cref{def:alg:j} terminates at the initial loop.
\end{proof}

Next, we study the asymptotic convergence of \cref{def:alg:j}.
For our problem, the lower-level value function $v(x)$
is always continuous under the {\it restricted inf-compactness} (RIC) condition
\cite[Definition 3.13]{GuoLinYeZhang}.
The RIC condition holds at $\bar{x}$ for $v(x)$
if $v(\bar{x})$ is finite and there exists a compact set $Z_0$ and
a scalar $\epsilon>0$ such that $S(x)\cap Z_0$ is nonempty for all
$x\in B_{\epsilon}(\bar{x})$ with $v(x)< v(\bar{x})+\epsilon$.
RIC conditions are very weak; see \cite{nie2017bilevel} for details.
For instance, RIC holds at $\bar{x}$ if $Z(x)$ is uniformly compact around $\bar{x}$.

\begin{lemma}
	For the lower-level optimization problem $(P_x)$,
	suppose RIC condition holds at $\bar x$ and
	$A$ is not the constant zero matrix.
	Then the value function $v(x)$ is continuous at $\bar x$.
\end{lemma}
\begin{proof}
Under the RIC condition, by \cite[Proposition 3.1]{GuoLinYeZhang23},
$v(x)$ is continuous at $\bar{x}$ if there exists $\bar y\in S(\bar x)$
such that
\[
	\lim_{x\rightarrow \bar x} \operatorname{dist} (\bar y, {Z}(x)) = 0,
\]
where $\operatorname{dist}$ stands for Euclidean distance.
Let $\bar y\in S(\bar x)$.
By Hoffman's Lemma and its proof in \cite{Hoffman},
there is a positive number $\kappa$ such that the error bound holds:
\begin{equation}\label{eq:distbd}
	\operatorname{{\rm dist}} (\bar y, {Z}(x))\leq \kappa \| (b(x) - A \bar y)_+\|,
\end{equation}
and the constant
\[
	\kappa \coloneqq
	\Big(\sum_{J\in\mc{P}} \|C(A_{J_0,J}\|^2 \Big)^{1/2}
	\Big(\sum_{J\in\mc{P}} (\det\,A_{J_0,J})^2  \Big)^{-1/2}.
\]
Here $J_0\in\mc{P}$ is a fixed label set,
$A_{J_0,J}$ is the submatrix with rows in $J_0$ and columns
in $J$, and $C(\cdot)$ denotes the cofactor matrix.
By \eqref{eq:distbd}, the $\operatorname{dist}(\bar y, {Z}(x))\to 0$ as $x\to \bar{x}$.
Hence, the continuity is proved.
\end{proof}

Let $\{ q^{(k)} \}_{k\in\N}$ with each $q^{(k)}: \re^n\times \re^p\rightarrow\re^p$
be a sequence of vector-valued continuous functions in $(x,y)$.
We say $\{ q^{(k)} \}_{k\in\N}$ is {\it equicontinuous}
at a given point $(\hat{x}, \hat{y})$ if for every $\rho>0$,
there exists $\epsilon>0$ such that
\begin{equation*}
	\| q^{(k)}(x,y) - q^{(k)}(\hat{x}, \hat{y})\|\le \epsilon,\quad
	\forall (x,y)\in B_{\rho}(\hat{x}, \hat{y})\,\,  \forall k\in\N.
\end{equation*}
When $\{q^{(k)}\}_{k\in\N}$ is a polynomial sequence,
a sufficient condition for the uniform continuity is that
$\{q^{(k)}\}_{k\in\N}$ has uniformly bounded degrees and coefficients.
\begin{theorem}\label{thm:asyconv:algj}
Under Assumption~\ref{ass:FE},
every accumulation point $(\hat{x}, \hat{y})$
of the sequence of $(x^{(k)}, y^{(k)})$ produced by
Algorithm~\ref{def:alg:j} is a global optimizer of $(P_J^*)$,
if $v(x)$ is continuous at $\hat{x}$ and $\{ q^{(k)}(x,y) \}_{k\in\N}$
is equicontinuous at $(\hat{x}, \hat{y})$.
\end{theorem}
\begin{proof}
For convenience, let $(P_{J,k})$ denote the optimization problem
\eqref{eq:Kj:alg} in the $k$-th loop.
Since $(\hat{x}, \hat{y})$ denote any accumulation point of $(x^{(k)}, y^{(k)})$,
we have $(\hat{x}, \hat{y})\in\mc{K}_J$ and up to selection of a subsequence,
we can generally assume $(x^{(k)}, y^{(k)}) \rightarrow (\hat{x}, \hat{y})$
as $k\rightarrow\infty$.
Since $F$ is a polynomial and each $(P_{J,k})$ is a relaxation of $(P_J^*)$,
we have
\[
	F(\hat{x}, \hat{y}) = \lim_{k\rightarrow\infty} F(x^{(k)}, y^{(k)})
	\le \vartheta_J^*.
\]
If $\hat{y}\in S(\hat{x})$, or equivalently, $v(\hat{x}) \ge f(\hat{x}, \hat{y})$,
then $(\hat{x}, \hat{y})$ is the global optimizer of $(P_J^*)$.
For every $k\in \N$, we can always decompose that
\[
	v(\hat{x}) - f(\hat{x}, \hat{y}) = S_1+S_2,
	\quad \text{where}
\]
\[
	S_1 \,\coloneqq\, v(\hat{x}) - f ( \hat{x}, q^{(k)}(\hat{x}, \hat{y}) ),\quad
	S_2 \,\coloneqq\, f( \hat{x}, q^{(k)}(\hat{x}, \hat{y}) ) - f(\hat{x}, \hat{y}).
\]
The point $(x^{(k)}, y^{(k)})$ must be feasible for all $(P_{J,s})$
with $s\le k$, so
\[
	f \big( x^{(k)}, q^{(s)}(x^{(k)}, y^{(k)}) \big)
	- f \big( x^{(k)}, y^{(k)} \big) \ge 0,
	\quad \forall s\in \N,\, s\le k.
\]
Let $k\rightarrow\infty$, the above implies $S_2\ge 0$,
thus $v(\hat{x}) - f(\hat{x}, \hat{y}) \ge S_1$.
Assume $v(x)$ is continuous at $\hat{x}$ and $\{q^{(k)}\}_{k\in \N}$ is
equicontinuous at $(\hat{x}, \hat{y})$.
By \cref{eq:fe:alg}, we have $v(x^{(k)}) = f ( x^{(k)}, q^{(k)} (x^{(k)}, y^{(k)}) )$.
As $k\rightarrow \infty$,
\[
	\begin{aligned}
	S_1 &= \big( v(\hat{x}) - v(x^{(k)}) \big) + \big( v(x^{(k)}) - f ( \hat{x}, q^{(k)}(\hat{x}, \hat{y}) \big)\\
	&=  \big( v(\hat{x}) - v(x^{(k)}) \big) + \big(f ( x^{(k)}, q^{(k)} (x^{(k)}, y^{(k)}) ) - f(\hat{x}, q^{(k)} (\hat{x}, \hat{y}) ) \big) \rightarrow 0,
\end{aligned}
\]
which implies $S_1+S_2\ge0$.
So $v(\hat{x})-f(\hat{x}, \hat{y})\ge 0$ and thus
$(\hat{x}, \hat{y})$ is a global optimizer of $(P_J^*)$.
\end{proof}

\section{Feasible extensions for linear constraints}
\label{sc:FE}

In this section, we introduce an efficient method to compute feasible extensions
for general linear or quadratic constraints.
That is, we seek a vector of polynomial vector
$q = (q_1,\ldots,q_p)$ satisfying \cref{eq:FE:Kj} for given
$(\hat{x}, \hat{y})\in\mc{K}_J\subseteq \mc{U}$ and $\hat{z}\in S(\hat{x})$.

\subsection{Linear feasible extensions}
First, we consider the case that
$b(x) = b_0 + Bx$ is a linear vector function.
Furthermore, we assume the upper level optimization
has linear constraints like
\[
	h(x,y) = h_0 + H_1x + H_2 y \ge 0,
\]
for a constant vector $h_0$ and two constant matrices $H_1, H_2$.
We look for the feasible extension function in the form
\[
	q(x, y) \, = \, w_0 + W_1x + W_2 y.
\]
Clearly, $q(x,y) \in Z(x)$ for all $(x,y) \in \mc{U}$
if it satisfies the polynomial identity
\[
	\baray{c}
	A( w_0 + W_1x + W_2 y) - b(x) \,=\, \xi + Y_1 h(x,y) + Y_2 (Ay - b(x) ) ,
	\earay
\]
for a nonnegative vector $\xi \ge 0$
and two nonnegative matrices $Y_1 \ge 0, Y_2 \ge 0$.
Thus, the feasible extension function $q$
can be obtained by solving the linear program:
\be \label{fe:lp}
	\left\{ \baray{c}
	w_0 + W_1\hat{x} +W_2\hat{y} = \hat{z},  \\
	A w_0 -b_0 = \xi + Y_1 h_0 -Y_2 b_0, \\
	A W_1 =  Y_1 H_1 - Y_2 B, \\
	A W_2 = Y_1 H_2 + Y_2 A, \\
	\xi \ge 0,\quad  Y_1 \ge 0,\quad  Y_2 \ge 0.
	\earay \right.
\ee

\begin{example}
Consider the BPOP with $x\in\re^2,y\in\re^2$ and
\[
	\begin{array}{l}
	h = (2-e^Tx,\, x_1,\, x_2),\\
	g = (-y_1+y_2+2x_1-2.5,\, -y_2-x_1+3x_2+2,\, y_1, y_2).
	\end{array}
\]
Let $\hat{x} = (1,1), \hat{y} = (3.5, 4), \hat{z} = (0, 4)$.
We look for a linear feasible extension by solving \cref{fe:lp}.
A feasible solution for \cref{fe:lp} is given by
\[
	w_0 = 0_{2\times 1},\quad
	W_1 = 0_{2\times 2},\quad
	\xi = 0_{4\times 1},\quad
	Y_1 = 0_{4\times 3},
\]
\[
	W_2 = \bbm 0 & 0 \\0 & 1\ebm, \quad
	Y_2 = \bbm 1 & 0 & 1 & 0\\
	 0 & 1 & 0 & 0\\
	 0 & 0 & 0 & 0\\
	 0 & 0 & 0 & 1\ebm.
\]
This gives the linear feasible extension $q(x,y) = (0, y_2)$.
\end{example}

\subsection{Quadratic feasible extensions}

Second, we consider the case that $h(x,y)$ and $b(x)$ are all quadratic polynomials.
We look for a quadratic feasible extension
$q(x, y) = (q_1(x,y), \ldots, q_p(x,y))$ such that
\[
	q_i(x, y) \, = \, (x \diamond y)^T  W_i (x \diamond y)
	\quad \text{where} \quad
	x \diamond y \coloneqq \bbm 1 & x^T & y^T \ebm^T .
\]
One can verify that $q(x,y) \in Z(x)$ for all $(x,y) \in \mc{U}$
if it satisfies the following system
(it is equivalent to a semidefinite program)
\be \label{fe:sdp}
	\left\{ \baray{rcl}
	q(\hat{x}, \hat{y} ) &=& \hat{z},  \\
	a_i^Tq(x,y) -b_i(x)  &=& (x \diamond y)^T  Y_i (x \diamond y)
	+  \nu_i^T h(x,y) + \theta_i^T (Ay - b(x)), \\
	\nu_i  &=& (\nu_{i,j} )_{ j \in \mc{E}_1 \cup \mc{I}_1 }, \,\,
	\nu_{i,j}  \ge  0 \, (j \in  \mc{I}_1 ), \\
	Y_i &\succeq & 0, \,\,  \theta_i \ge 0, \, i =1, \ldots, p.
	\earay  \right.
\ee
We remark that the second equation in \cref{fe:sdp} means that
the left and right hand side polynomials in $(x,y)$ are the same.

\begin{example}
Consider the BPOP with $x_1\in\re^2, y\in\re^2$ and
\[
	\begin{array}{l}
	h = (4-x_1^2-2x_2,\, x_1,\, x_2),\\
	g = (-2y_1+y_2+x_1^2-2x_1+x_2^2+3,\,
	3y_1-4y_2+x_2-4,\, y_1,\, y_2).
	\end{array}
\]
Let $\hat{x} = (1,1), \hat{y} = (1.5, 0), \hat{z} = (1.8,0.6)$.
We look for a quadratic feasible extension by solving \cref{fe:sdp}.
A feasible solution for \cref{fe:sdp} is given by
\[
	\begin{array}{c}
	W_1 = \left[\begin{array}{rrccc}
	1.6 & -0.8 & 0.1 & 0 & 0\\
	-0.8 & 0.8 & 0 & 0 & 0\\
	0.1 & 0 & 0.8 & 0 & 0\\
	0 & 0 & 0 & 0 & 0\\
	0 & 0 & 0 & 0 & 0
	\end{array}\right],\quad
	W_2 = \left[\begin{array}{rrccc}
	0.2 & -0.6 & 0.2 & 0 & 0\\
	-0.6 & 0.6 & 0 & 0 & 0\\
	0.2 & 0 & 0.6 & 0 & 0\\
	0 & 0 & 0 & 0 & 0\\
	0 & 0 & 0 & 0 & 0
	\end{array}\right],\\
	\theta_1 = \theta_2 = \bbm 0\\0\\0\\0\ebm,\,
	\theta_3 = \bbm 0.8\\0.2\\1\\0\ebm,\,
	\theta_4 = \bbm 0.6\\0.4\\0\\1\ebm,
	Y_1 = Y_2 = 0_{5\times 5}, \, \nu_1 = \nu_2 = 0_{3\times 1} .
	\end{array}
\]
This gives the following quadratic feasible extension function $q(x,y)=$
\[
	(0.8x_1^2-1.6x_1+0.2x_2+0.8x_2^2+1.6,\,
	0.6x_1^2-1.2x_1+0.4x_2+0.6x_2^2+0.2).
\]
\end{example}

For more general constraints, it is typically a challenge
to get feasible extensions. The problem is mostly open,
to the best of our knowledge.
This is an interesting question for future work.

\section{Numerical examples}\label{sc:numexam}

In this section, we give numerical experiments to solve BPOPs with
\cref{def:alg:bpop}.
In computations, each optimization problem is solved globally by
Moment-SOS relaxations with software {\tt GloptiPoly 3} \cite{GloPol3},
\texttt{Mosek} \cite{mosek}
and {\tt SeDuMi} \cite{sturm1999using}.
We implement algorithms using \texttt{MATLAB R2024a},
in a laptop with CPU 8th Generation Intel® Core™ Ultra 9 185H and RAM 32 GB.
In the following, we report all numerical results with four digits.
In each problem, the constraints are ordered from left to
right, and from top to bottom.
The CPU time is given with the unit ``second''.
Suppose $(x^*, y^*)$ is the computed solution for \eqref{eq:genblv}
by some algorithm (we compared \cref{def:alg:bpop} with some other methods).
Its active set $\mc{I}(x^*,y^*)$ is computed
within a small tolerance, e.g.,
$\{i\in\mc{I}\,\vert\, a_i^Tx^*-b_i(y^*)\le 10^{-4} \}$.

\subsection{Some explicit examples}
We present some explicit examples of BPOPs.

\begin{example}\label{ex:cons_grag}
Consider the BPOP
\[
	\left\{\begin{array}{cl}
	\min\limits_{x\in\re^2,y\in\re^3}
	& (x_1-x_2)^2y_3-2x_2^2y_2y_3-y_1^2y_2y_3\\
	\st & 4-\|x\|^2-\|y\|^2\ge 0,\\
	& \|x\|^2 - 1\ge 0,\,x \ge 0\\
	&  y\in S(x),
	\end{array}\right.
	\]
where $S(x)$ is the set of global optimizer(s) of
\[
	\left\{\begin{array}{cl}
	\min\limits_{y\in\re^3}
	& (y_1+2y_2-x_1y_3)^2+x_1y_1+x_2y_2\\
	\st & y\ge 0,\,\, y_1 - y_2 - y_3 + 1 \ge 0,\\
	& y_1 - 2y_2 + 0.5y_3 - (2x_1-1) \ge 0,\\
	& -2y_1 + y_2 + 0.5y_3 -(2x_2-1) \ge 0.\\
	\end{array}\right.
\]
Since the lower-level problem is convex,
we do not need to compute feasible extensions.
It took 1.68 second for \cref{def:alg:bpop} to solve this problem.
There are total $|\mc{P}| = 20$ branch problems.
The index set for feasible branch problems is
\[
	\begin{array}{c}
	\mc{J}_0 =  \{ \{1,2,4\}, \{1,5,6\}, \{2,3,5\},
	\{2,4,5\},\{2,5,6\}, \{3,5,6\}, \{4,5,6\}\}.
	\end{array}
\]
We got the global optimal value and optimizer
\[
	F_{min} = -0.0876,\quad
	x^* = (0.6852, 0.7284),\,
	y^* = (0.0576, 0.0864, 0.9712),
\]
from the branch problem with index contained in
$\mc{J}_1 = \{\{4,5,6\}\}$.
We did not find non-global local optimizer,
i.e., $\mc{J}_2 = \emptyset$.
\end{example}

\begin{example}\label{ex:LnoncMatcon}
Consider the BPOP
\begin{equation}\label{eq:lnonc}
	\left\{\begin{array}{cl}
		\min\limits_{x\in \re^4,y\in\re^2} & \|x-e\|^2 + \|y+e\|^2\\
		\st & x\ge 0,\quad y_1+1\ge 0,\quad y_2\ge 0,\\
		& 4-e^Tx-e^Ty\ge 0,\quad  y\in S(x),
	\end{array}\right.
\end{equation}
where $S(x)$ is the set of global optimizer(s) of
\[
	\left\{\begin{array}{cl}
	\min\limits_{y\in\re^2}
	& (x_3-y_1-1)(x_4-y_2+1)-\|y\|^2,\\
	\st & y_1 -x_1 +1 \ge 0,\quad y_1+y_2+x_1-1\ge 0,\\
	& 1.5-y_1-x_1\ge 0,\quad\,\,\, -1+x_2+y_2\ge 0,\\
	& 3-x_1-y_1-y_2\ge 0.
	\end{array}\right.
\]
Since the lower-level problem of \eqref{eq:lnonc} is nonconvex,
we need to use feasible extensions in computation.
We listed feasible extensions for computed tuples $(\hat{x}, \hat{y})\in \mc{K}$
and $\hat{z}\in S(\hat{x})$ in \cref{tab:fe}.
\begin{table}[htb]
	\small
	\centering
	\caption{Feasible extensions used for \cref{ex:LnoncMatcon}}
	\label{tab:fe}
	\begin{tabular}{|c|c|c|}
		\hline
		$(\hat{x}, \hat{y})$ & $\hat{z}$ & $q(x,y)$\\
		\hline
		$(1.2500, 1.0000, 0.7500, 0.7500),(0.2500, 0.0000)$
		& $(0.2500, 1.5000)$ & $(0.25, 2.75-x_1)$\\
		\hline
		$((0.9337, 1.0000, 0.0000, 1.4356), (0.5663, 0.0000)$
		& $(-0.0663, 2.1326)$ & $(x_1-1, 4-2x_1)$\\
		\hline
	\end{tabular}
\end{table}
One can easily verify that these $q$ satisfies \cref{ass:FE} since
\[
	\begin{array}{rcl}
	1.5-x_1\ge&  y_1 & \ge x_1-1,\\
	3-x_1-y_1\ge &  y_2 & \ge \max\{1-x_1-y_1,\, x_2-1\},
	\end{array}
\]
for every feasible point of \eqref{eq:lnonc}.
It took 1.68 second for \cref{def:alg:bpop} to solve \eqref{eq:lnonc}.
There are total $|\mc{P}| = 8$ branch problems and only one of them is feasible,
i.e., $\mc{J}_0 = \mc{J}_1 = \{\{1,5\}\}$.
We got the global optimal value and optimizer.
\[
	\begin{array}{ll}
	F_{\min} = 9.2083,& x^* = (1.2500,\, 0.3336,\, 0.3332,\, 0.3332),\\
	& y^* = (0.2500,\, 1.5000).
	\end{array}
\]
\end{example}

\begin{example}\label{ex:lin_simplex}
Consider the BPOP
\begin{equation}\label{eq:lin_simplex}
	\left\{\begin{array}{cl}
		\min\limits_{x,y\in\re^4}
		& x_1^2(y_3^2-1)+x_4^2(y_1^2-2)-x_2x_4+x_1x_3\\
		\st & 1\le e^Tx\le 4,\,  2x_1x_2-y_1^2-2y_2^2\ge 0,\\
		& x_1 - y_1\ge 0,\,  x_2 - y_2\ge 0,\, 2-x_1-x_2\ge 0,\\
		& 2-x_3^2-x_4^2\ge 0,\quad  y\in S(x),
	\end{array}\right.
\end{equation}
where $S(x)$ is the set of global optimizer(s) of
\[
	\left\{\begin{array}{cl}
	\min\limits_{y\in\re^4}
	& x_1y_1^2 - x_2y_2^2 + x_3(y_3-1)^2 + x_4(y_4+1)^2\\
	\st & y_1- y_2 - x_2 \ge 0,\, x_1 - y_1 + y_2 \ge 0,\\
	& y_1 + y_2 + x_1 +x_2 \ge 0,\\
	& 4x_1 - 2x_2-y_1 - y_2\ge 0,\\
	& y_3\ge 0,\quad y_4\ge 0,\quad 3-y_3 -y_4 \ge 0.
	\end{array}\right.
\]
Since the lower-level problem of \eqref{eq:lin_simplex} is nonconvex,
we may need feasible extensions in computations.
For $(\hat{x},\hat{y})\in\mc{K}$ and $\hat{z}\in S(\hat{x})$,
one can verify that
\[
	q = \frac{1}{2}\bbm (\mu_1+5\mu_2-1)x_1-(\mu_1+\mu_2)x_2\\
	(5\mu_2-\mu_1-1)x_1+(\mu_1-5\mu_2-2)x_2\\
	2\hat{z}_3\\ 2\hat{z}_4
	\ebm\quad \mbox{with}\quad
	\begin{aligned}
	&	\mu_1 = \frac{\hat{z}_1-\hat{z}_2-\hat{x}_2}{\hat{x}_1-\hat{x_2}},\\
	&	\mu_2 = \frac{\hat{z}_1+\hat{z}_2+\hat{x}_1+\hat{x}_2}{5\hat{x}_1-\hat{x}_2},
	\end{aligned}
\]
satisfies Assumption~\ref{ass:FE}.
It took $2.65$ seconds for \cref{def:alg:bpop} to solve this problem.
There are total $|\mc{P}| = 12$ branch problems for \eqref{eq:lin_simplex}.
All of them are feasible, i.e., $\mc{J}_0 = \mc{P}$.
We got the global optimal value and optimizer
\[
	\begin{array}{ll}
	F_{min} = -8.0000,
	&  x^* = (2.0000, 0.0000, 0.0000, 1.4142),\\
	& y^* = (1.4142, 0.0000, 0.0000, 0.0000),
	\end{array}
\]
from branch problems with index contained in
\[
	\begin{array}{c}
	\mc{J}_1 = \{\{1,3,5,6\},\, \{1,3,6,7\}, \, \{1,4,5,6\},\, \{1,4,6,7\},\\
	\{2,3,5,6\},\,\{2,3,6,7\},\,\{2,4,5,6\},\, \{2,4,6,7\}\}.
	\end{array}
\]
In addition, we got a non-global local optimal value and optimizer
\[
	\begin{array}{ll}
	\vartheta_J^* = -7.0000,
	& u_J^* = (2.0000, 0.0000, 0.3660, -1.3660),\\
	& v_J^* = (0.0000, 0.0000, 0.0000, 3.0000),
	\end{array}
\]
from branch problems with index contained in
\[
	\begin{array}{c}
	\mc{J}_2 = \{ \{1,3,5,7\},\, \{1,4,5,7\},\, \{2,3,5,7\},\, \{2,4,5,7\}\}.
	\end{array}
\]
\end{example}

\begin{example}\label{ex:toll}
The toll-setting problem in transportation can be formulated into a bilevel
optimization problem.
Consider the toll-setting bilevel problem from \cite{Kue2023}:
\[
	\left\{\begin{array}{cl}
	\min\limits_{x,y\in\re^8} & (5-x_3)y_3-x_4y_4-x_8y_8\\
	\st & x_1 = 2,\, x_2 = 6,x_3\ge 5,\, x_4\ge 0,\\
	& x_5=4,\, x_6 = 2,\, x_7=6,\,x_8\ge 0,\\
	& y\in S(x),
	\end{array}\right.
\]
where $S(x)$ is the optimizer set of
\[
	\left\{\begin{array}{cl}
	\min\limits_{y\in\re^8} & x^Ty\\
	\st & y_1+y_2+y_3 = 1,\, y_8+y_7+y_3 = 1,\\
	& y_4+y_5-y_1 = 0,\, y_6+y_7-y_2-y_4 = 0,\\
	& y_8-y_5-y_6 = 0,\, y_i\ge 0,\,\forall i\in[8].
	\end{array}\right.
\]
Since the lower-level problem is convex, we do not need to compute feasible extensions.
We used the equality constraints to eliminate variables and reformulate the BPOP.
It took 10.95 seconds for \Cref{def:alg:bpop} to solve total $|\mc{P}| = 45$
branch problems for the reformulated BPOP,
where $|\mc{J}_0| = 36$, $|\mc{J}_1| = 1$ and $|\mc{J}_2| = 26$.
We obtain the global optimal value and the global optimizer
\[
	\begin{array}{l}
	F_{min} = -7.0176, \quad
	x^* = (2,6,8.0595,4.0000, 4,2,6,4.0000),    \\
	y^* = (1.0000, 0.0000, 1.0000,0.0000,
	1.0000, 0.0000, 0.0000, 1.0000).
	\end{array}\]
\end{example}

\subsection{Comparison with other methods}

We compare \cref{def:alg:bpop} with some existing methods for solving BPOPs
such as the LME-FE method \cite{nie2020bilevel} and the KKT approach.
For the KKT approach, we use \texttt{BARON} \cite{Baron} to solve $(P_{kkt})$ with
the starting point $(x,y,\lambda)$ where $\lambda=0$ and $(x,y)$ is
randomly chosen.	
We allocated a maximum CPU time allowance of 500 seconds for each method per problem.
For a computed point $(x^*,y^*)$, we evaluate its lower-level feasibility by
\[
	\eta^*\,\coloneqq\,f(x^*,y^*)-v(x^*),
\]
where $v(x)$ is the lower-level value function defined in \eqref{eq:valfun}.

We applied LME-FE method and the KKT approach to Examples~\ref{ex:cons_grag}--\ref{ex:toll}
and reported numerical results in \cref{tab:compare}.
From the table,
we can see that the LME-FE method failed to converge for Example~\ref{ex:toll},
when the maximum time allowance is reached.
The corresponding moment problem is very difficult to solve since
the LME implied by \cref{eq:Lnd} has a high degree (e.g., $\deg(H(x,y) = 16$).
The KKT approach falsely returns ``infeasible'' for Example~\ref{ex:LnoncMatcon}
and failed to converge to a feasible point for
Example~\ref{ex:lin_simplex} and Example~\ref{ex:toll}.

\begin{table}[htb!]
\small
\centering
\caption{Comparison with some existing methods}
\label{tab:compare}
\begin{tabular}{|c|l|l|c|c|}
	\hline
	Method & $F^*$ & $(x^*,y^*)$ & $\eta^*$ & Time \\
	\hline
	\multicolumn{5}{|c|}{{\bf \Cref{ex:cons_grag}}}\\
	\hline
	Alg.~\ref{def:alg:bpop} & $-0.0876$ & $\begin{array}{l}
		(0.6852,0.7284),\\
		(0.0576,0.0864,0.9712)\end{array}$ & $-3.55\cdot10^{-9}$ & 1.68\\
	\hline
	LME-FE & $-0.0876$ & $\begin{array}{l}
		(0.6852, 0.7284)\\
		(0.0576, 0.0864, 0.9712)
	\end{array}$ & $-7.77\cdot10^{-8}$ & 8.51\\
	\hline
	$(P_{kkt})$ & $-0.0876$ & $\begin{array}{l}
		(0.6852,0.7284),\\
		(0.0577,0.0864,0.9712)
	\end{array}$ & $-3.23\cdot10^{-10}$ & 500\\
	\hline
	\multicolumn{5}{|c|}{{\bf \Cref{ex:LnoncMatcon}}}\\
	\hline
	Alg.~\ref{def:alg:bpop} & $9.2083$ & $\begin{array}{l}
		(1.2500, 0.3336, 0.3332, 0.3332),\\
		(0.2500, 1.5000)
	\end{array}$ & $-1.00\cdot10^{-6}$ & 1.68\\
	\hline
	LME-FE & $9.2083$ & $\begin{array}{l}
		(1.2500, 0.3332, 0.3334, 0.3334),\\
		(0.2500, 1.5000)
	\end{array}$ & $-1.72\cdot 10^{-9}$ & 13.31\\
	\hline
	$(P_{kkt})$ & \multicolumn{4}{c|}{not convergent}\\
	\hline
	\multicolumn{5}{|c|}{{\bf \Cref{ex:lin_simplex}}}\\
	\hline
	Alg.~\ref{def:alg:bpop} & $-8.0000$ & $\begin{array}{l}
		(2.0000, 0.0000, 0.0000, 1.4142),\\
		(0.0000, 0.0000, 0.0000, 0.0000)
	\end{array}$ & $-5.07\cdot10^{-7}$ & 3.62\\
	\hline
	LME-FE & $-8.0000$ & $\begin{array}{l}
		(2.0000, 0.0000, 0.0000, 1.4142),\\
		(0.0000, 0.0000, 0.0000, 0.0000)
	\end{array}$ & $-7.41\cdot 10^{-8}$ & 4.53\\
	\hline
	$(P_{kkt})$ & $-8.5000$ & $\begin{array}{l}
		(2.0000, 0.0000, -0.5000, 1.3229),\\
		(0.0000, 0.0000, 1.0000, 0.0000)
	\end{array}$ & $-2.00$ & 0.14\\
	\hline
	\multicolumn{5}{|c|}{{\bf \Cref{ex:toll}}}\\
	\hline
	Alg.~\ref{def:alg:bpop} & $-7.0176$ & $\begin{array}{l}
		(2,6,11.9667,4.3783,
		4,2,6,8.0503),\\
		(0.0557, 0.0550,0.0000, 0.9441,\\
		0.0546,0.0000,0.0007,0.0552)
	\end{array}$ & $-2.65\cdot10^{-8}$ & 10.95\\
	\hline
	LME-FE & \multicolumn{4}{c|}{not convergent}\\
	\hline
	$(P_{kkt})$ & $-8.0000$ & $\begin{array}{l}
		(2,6,8.0595,4.0000,
		4,2,6,4.0000),\\
		(1.0000, 0.0000, 1.0000,0.0000,\\
		1.0000, 0.0000, 0.0000, 1.0000)
	\end{array}$ & $-3.9405$ & 0.06\\
	\hline
\end{tabular}
\end{table}

In addition, we compare \cref{def:alg:bpop}
with the LME-FE method and the KKT approach for some randomly generated BPOPs.

\begin{example}\label{ex:Rplus}
Consider BPOPs with ball constraints on $x$ and
nonnegative constraints on $y$:
\begin{equation}\label{eq:Rplus}
	\left\{\begin{array}{cl}
		\min\limits_{x\in\re^n, y\in\re^p} & c_1^T\bbm x\\y\ebm_{3}+
		\left\|B_1\bbm x\\y\ebm^{2}\right\|^2\\
		\st & \|x\|^2\le 1,\quad y\in S(x),
	\end{array}\right.
\end{equation}
where $S(x)$ is the optimizer set of
\[
	\left\{\begin{array}{cl}
	\min\limits_{y\in\re^p} & c_2^T[x]_{3}+c_3^T[y]_{3}+
	\left\|B_2 \bbm [x]^{2}\\ [y]^{2}\ebm \right\|^2\\
	\st & y\ge 0.
	\end{array}\right.
\]
In the above, $[x]_d$ is the vector of monomials in $x$ with degrees up to $d$
and $x^{[d]}$ is the vector of monomials in $x$ with degrees equal to $d$.
The symbol is defined similarly when $x$ is changed to other variables such as
$y$ or $(x,y)$.
The \cref{eq:Rplus} is a simple bilevel problem.
For computed points $(\hat{x}, \hat{y})\in\mc{K}$ and $\hat{z}\in S(\hat{x})$,
we can choose $q = \hat{z}$ for \cref{ass:FE} to satisfy.
Since \eqref{eq:Rplus} has the same dimension for the lower-level constraints
and variables, its disjunctive program reformulation equals $(P)$ and
\cref{def:alg:bpop} reduces to the LME-FE approach.
The computational results are reported in \Cref{tab:Rplus}.
In the table, we randomly generated $10$ instances for each case.
These BPOPs are successfully solved by \Cref{def:alg:bpop}.
On the other hand, \texttt{BARON} did not converge to a feasible point of
BPOPs at any instance,
when the maximum CPU time allowance of 500 seconds was reached.
It shows that our \cref{def:alg:bpop} is much more efficient than
the KKT approach.
\begin{table}[htb]
	\small
	\centering
	\caption{Numerical results of \Cref{ex:Rplus}}
	\label{tab:Rplus}
	\begin{tabular}{|c|c|c|c|c|c|c|c|}
		\hline
		\multirow{2}{*}{Method} & \multirow{2}{*}{$(n,p)$} &
		\multicolumn{3}{|c|}{Time} & \multicolumn{3}{|c|}{$\eta^*$}\\
		\cline{3-8}
		& & \texttt{Min} & \texttt{Avg} &
		\texttt{Max} & \texttt{Min} & \texttt{Avg} & \texttt{Max}\\
		\hline
		\multirow{3}{*}{Alg~\ref{def:alg:bpop}} & $(3,3)$  & 0.23 & 0.91 & 3.84 & $-7.06\cdot 10^{-6}$ & $-1.03\cdot 10^{-6}$ & $-1.48\cdot10^{-8}$\\
		\cline{2-8}
		& $(4,4)$ & $1.46$ & $6.58$ & $39.43$ & $-5.90\cdot10^{-6}$ & $-8.07\cdot10^{-7}$
		& $-2.05\cdot10^{-8}$\\
		\cline{2-8}
		& $(5,5)$  & $3.79$ & $14.48$ & $54.68$ & $-9.41\cdot 10^{-8}$ & $-3.83\cdot10^{-8}$ & $-1.25\cdot10^{-8}$ \\
		\hline
		\multirow{2}{*}{$(P_{kkt})$} & $(3,3)$  & 500 & 500 & 500 & $-0.5287$ & $-0.2771$ & $-0.0067$\\ \cline{2-8}
		& $(4,4)$ & 500 & 500 & 500 & $-0.1571$ & $-0.0781$ & $-0.0324$\\
		\hline
	\end{tabular}
\end{table}
\end{example}

\subsection{Some examples from references}
We use \Cref{def:alg:bpop} to solve all test problems in
the bilevel optimization library \cite{BOLIB} that has
linear lower-level constraints.
All problems are solved accurately up to a numerical error,
i.e., $\eta^*$.
We list numerical results in \Cref{tab:BOLIB},
where $F^*$ denotes the computed optimal value and
$(x^*,y^*)$ denotes the computed optimizer.
The $|\mc{P}|$ denotes the total number of branch problems,
$|\mc{J}_0|$ denotes the number of feasible branch problems
and $|\mc{J}_1\cup \mc{J}_2|$ denotes the number of branch problems whose optimizers
are verified to be global/local optimizers of the original bilevel problem.
For BPOPs with parameters,
we choose $\rho=3$ for the problem \texttt{CalamaiVicente1994a},
$c=1$ for the problem \texttt{HenrionSurowiec2011},
and $M=2$ for the problem \texttt{IshizukaAiyoshi1992a}.

{\scriptsize
\begin{longtable}{|l|c|c|l|c|c|}
	\caption{Numerical results for test problems in bilevel optimization library}
	\label{tab:BOLIB}\\
	\hline
	{\bf Problem name} & $F_{min}$ & $(x^*,y^*)$ & $|\eta^*|$ &
	$\begin{array}{c}(|\mc{P}|, |\mc{J}_0|,\\ |\mc{J}_1\cup \mc{J}_2|)\end{array}$ & Time \\
	\hline
	\endhead	
	\hline
	\endfoot
	AiyoshiShimizu & $0.0000$  &  $
	(0.0000,30.0000),\,(-10.0000,10.0000)$ & $1.79\cdot10^{-6}$
	& (9,9,9) & 0.46 \\
	\hline
	AllendeStill2013 & $1.0000$ & $(0.5000, 0.5000),\,
	(0.5000, 0.50000)$ & $1.22\cdot10^{-8}$ & (4,4,4) & 0.09\\
	\hline
	AnEtal2009 & $2.2516\cdot10^3$ & $(0.2000, 2.0000),\,
	(4.0000, 4.6000)$ & $3.30\cdot10^{-5}$ & (6,1,1) & 0.28\\
	\hline
	Bard1988Ex1 & $17.0000$ & $(1.0000),\, (0.0000)$ & $9.38\cdot10^{-9}$ & (4,4,2) & 0.17\\
	\hline
	Bard1988Ex2 & $-6600$ & $\begin{array}{c}(7.0000, 3.0000,12.0000,18.0000),\\
		(0.0000, 10.0000, 30.0000,0.0000)\end{array}$ & $1.26\cdot 10^{-6}$ &
	(169,42,22) & 13.51\\
	\hline
	Bard1988Ex3 & $-12.6787$ & $\begin{array}{c}(0.0005,2.0000),\\
		(1.8750,0.9062)\end{array}$ & $1.26\cdot 10^{-6}$ &
	(6,3,3) & 0.04\\
	\hline
	Bard1991Ex1 & $2.0000$ & $(2.0000),\,
	(6.0000, 0.0000)$ & $1.37\cdot10^{-7}$ &
	(3,2,2) & 0.11\\
	\hline
	BardBook1998 & $0.0000$ & $\begin{array}{c}(30.0000, 40.0000),\\
		(10.0000, 20.0000)\end{array}$ & $3.62\cdot10^{-7}$ &
	(15,4,2) & 0.45\\
	\hline
	CalamaiVicente1994a & $0.2500$ & $(0.5000),\,(0.5000)$ & $2.44\cdot10^{-9}$ &
	(3,3,3) & 0.03\\
	\hline
	CalamaiVicente1994b & $0.3125$ & $\begin{array}{c}
		(1.2500, 0.50001.0000, 1.0000)\\
		(0.2500, 0.5000)\end{array}$ & $1.22\cdot10^{-8}$ &
	(9,9,9) & 0.22\\
	\hline
	CalamaiVicente1994c & $0.3125$ & $\begin{array}{c}
		(0.1308,0.0520,0.1022,0.0674)\\
		(0.0250, 0.0250)
	\end{array}$ & $1.36\cdot 10^{-7}$ &
	(9,9,9) & 0.34\\
	\hline
	ClarkWesterberg1990a & $5.0000$ & $(1.0000),\,(3.0000)$ & $1.69\cdot10^{-6}$ &
	(3,3,3) & 0.09\\
	\hline
	Colson2002BIPA1 & $250.0000$ & $(5.0000),\, (5.0000)$ & $1.92\cdot10^{-7}$ &
	(3,3,3) & 0.09\\
	\hline
	Colson2002BIPA2 & $17.0000$ & $(1.0000),\, (0.0000)$ & $3.13\cdot10^{-9}$ &
	(4,4,2) & 0.14\\
	\hline
	Colson2002BIPA4 & $0.0469$ & $(0.0001),\, (0.5774)$ & $1.88\cdot10^{-9}$ & (2,0,2) & 0.10 \\
	\hline
	DempeEtal2012 & $-1.0000$ & $(-1.0000),\, (1.0000)$ & $1.75\cdot 10^{-9}$ &
	(2,2,2) & 0.12\\
	\hline
	DempeFranke2001Ex41 & $5.0000$ & $(0.0000, -1.0000),\,(1.0000, 2.0000)$ & $1.77\cdot 10^{-7}$ &
	(5,2,2) & 0.32\\
	\hline
	DempeFranke2014Ex42 & $2.1250$ & $(-1.0000, -1.0000),\, (2.2500,1.2500)$ & $1.87\cdot10^{-7}$ &
	(2,2,2) & 0.10\\
	\hline
	DempeFrank2014Ex38 & $-1.0000$ & $(0.0000, -1.0000),\, (1.0000, 2.0000)$ & $1.47\cdot 10^{-7}$ &
	(3,2,2) & 0.31\\
	\hline
	DempeLohse2011Ex31a & $-5.7500$ & $(0.0000, 0.5000),\, (2.0000, 0.0000)$ & $5.55\cdot 10^{-9}$ &
	(6,6,6) & 0.29\\
	\hline
	DempeLohse2011Ex31b & $-12.0000$ & $\begin{array}{l}(0.5000, 0.5000,0.0000)\\
		(0.0000, 0.0000, 2.0000)\end{array}$ & $2.08\cdot10^{-8}$ &
	(9,9,4) & 0.85\\
	\hline
	DeSilva1978 & $-1.0000$ & $(0.5000, 0.5000),\, (0.5000, 0.5000)$ & $2.97\cdot10^{-9}$ &
	(4,4,4) & 0.48\\
	\hline
	FalkLiu1995 & $-2.2500$ & $(0.7500, 0.7500),\, (0.7500, 0.7500)$ & $5.20\cdot10^{-8}$ &
	(4,4,4) & 0.18\\
	\hline
	FloudasEtal2013 & $0.0000$ & $(0.0000, 30.0000),\, (-10.0000, 10.0000)$ & $9.53\cdot10^{-8}$ &
	(15,15,12) & 0.64\\
	\hline
	GumusFloudas2001Ex1 & $2.250\cdot10^3$ & $(11.2500),\, (5.0000)$ & $-6.76\cdot10^{-8}$ &
	(3,3,3) & 0.14\\
	\hline
	GumusFloudas2001Ex4 & $9.0000$ & $(3.0000),\, (5.0000)$ & $4.61\cdot10^{-8}$ &
	(2,2,2) & 0.10\\
	\hline
	HatzEtal2013 & $1.4897\cdot10^{-9}$ & $(0.0000),\, (0.0000, 0.0000)$ & $2.39\cdot 10^{6}$ &
	(1,1,1) & 0.04\\
	\hline
	HendersonQuandt1958 & $-3.1681\cdot 10^3$ & $(1.0951),\, (22.6228)$ & $0.0593$ &
	(1,1,1) & 0.01\\
	\hline
	HenrionSurowiec2011 & $-0.0625$ & $(-0.2500),\, (-0.2500)$ & $6.00\cdot10^{-9}$ & (1,1,1) & 0.09 \\
	\hline
	IshizukaAiyoshi1992a & $1.7848\cdot10^{-9}$ & $0.1359,\, (-0.1359,0.0023)$ & $4.47\cdot10^{-10}$ &
	(4,2,2) & 0.50\\
	\hline
	KleniatiAdjiman2014Ex3 & $-1.0000$ & $(0.0000),\, (1.0000)$ & $6.02\cdot10^{-9}$ & (2,2,2) & 0.06 \\
	\hline
	LamparielloSagratellaEx31 & $1.0000$ & $(1.0000),\, (0.0000)$ & $7.59\cdot 10^{-9}$ & (1,1,1) & 0.09\\
	\hline
	LamparielloSagratellaEx32 & $0.5000$ & $(0.5000),\,(0.5000)$ & $2.93\cdot 10^{-9}$ & (1,1,1) & 0.01 \\
	\hline
	LamparielloSagratellaEx33 & $0.5000$ & $(0.5000),\, (0.0000, 0.5000)$ & $3.09\cdot 10^{-9}$ &
	(3,2,1) & 0.29\\
	\hline
	LamparielloSagratellaEx35 & $0.8000$ & $(0.8000),\, (0.4000)$ & $1.09\cdot 10^{-9}$ & (3,2,2) & 0.15\\
	\hline
	LucchettiEtal1987 & $4.6278\cdot10^{-10}$ & $(1.0000),\, (0.0000)$ & $9.50\cdot 10^{-10}$ & (2,2,2) & 0.09\\
	\hline
	LuDebSinha2016f & $-8.1522$ & $(105.2174, 8.1522),\, (2.5000)$ & $6.02\cdot 10^{-9}$ & (1,1,1) & 0.04 \\
	\hline
	MacalHurter1997 & $81.3279$ & $(10.0164),\, (0.8201)$ & $6.64\cdot10^{-9}$ & (1,1,1) & 0.01 \\
	\hline
	MitsosBarton2006Ex39 & $-1.0000$ & $(-1.0000),\, (-1.0000)$ & $1.14\cdot 10^{-9}$ & (2,1,1) & 0.03 \\
	\hline
	MitsosBarton2006Ex310 & $0.5000$ & $(0.4380),\,(0.5000)$ & $2.64\cdot 10^{-9}$ & (2,2,2) & 0.25 \\
	\hline
	MitsosBarton2006Ex311 & $-0.8000$ & $(0.0000),\,(-0.8000)$ & $1.40\cdot 10^{-8}$  & (2,2,2) & 0.15 \\
	\hline
	MitsosBarton2006Ex312 & $-0.0000$ & $(0.0000),\, (0.0000)$ & $2.97\cdot 10^{-9}$ & (2,2,2) & 0.03 \\
	\hline
	MitsosBarton2006Ex313 &  $-1.0000$ & $(0.0000),\, (-1.0000)$ & $2.29\cdot 10^{-10}$ & (2,2,1) & 0.21 \\
	\hline
	MitsosBarton2006Ex314 & $0.2500$ & $(0.2500),\, (0.5000)$ & $5.25\cdot 10^{-11}$ & (2,2,2) & 0.23 \\
	\hline
	MitsosBarton2006Ex315 & $0.0000$ & $(-1.0000),\,(1.0000)$ & $6.40\cdot 10^{-9}$  & (2,2,2) & 0.23 \\
	\hline
	MitsosBarton2006Ex316 & $-2.0000$ & $(-1.0000),\, (0.0000)$ & $1.05\cdot 10^{-9}$ & (2,2,2) & 0.18 \\
	\hline
	MitsosBarton2006Ex317 & $0.1542$  & $(-0.2373),\, (0.0000)$ & $9.53\cdot 10^{-5}$ & (2,2,2) & 0.53 \\
	\hline
	MitsosBarton2006Ex318 & $-1.0000$ & $(1.0000),\, (0.0000)$ & $1.27\cdot 10^{-9}$  & (2,2,2) & 0.28 \\
	\hline
	MitsosBarton2006Ex319 & $-0.2581$ & $(0.1886),\, (0.4343)$ & $1.23\cdot 10^{-9}$  & (2,2,2) & 0.31 \\
	\hline
	MitsosBarton2006Ex320 & $0.3125$ & $(0.5000),\,(0.5000)$ & $3.20\cdot10^{-10}$ & (2,2,2) & 0.23 \\
	\hline
	MitsosBarton2006Ex321 & $0.2095$ & $(-0.5545),\, (0.4554)$ & $1.38\cdot 10^{-7}$  & (2,2,2) & 0.18 \\
	\hline
	MitsosBarton2006Ex324 & $-1.7547$ & $(0.2106),\, (1.7991)$ & $2.62\cdot 10^{-7}$ & (2,2,2) & 0.09 \\
	\hline
	MitsosBarton2006Ex326 & $-2.3536$ & $\begin{array}{l}(-1.0000, -1.0000),\\
		(1.0000, -1.0000, -0.7071)\end{array}$ & $3.90\cdot 10^{-8}$ & (8,8,8) & 6.09 \\
	\hline
	MorganPatrone2006a & $-1.0000$ & $0.0000,\, 1.0000$ & $1.94\cdot 10^{-10}$ & (2,2,2) & 0.09\\
	\hline
	MorganPatrone2006b & $-1.2500$ & $0.2500,\, 1.0000$ & $3.84\cdot 10^{-10}$ & (2,2,2) & 0.06 \\
	\hline
	MorganPatrone2006c & $-1.0000$ & $2.0000,\, -1.0000$ & $1.94\cdot 10^{-8}$ & (2,1,1) & 0.03 \\
	\hline
	MuuQuy2003Ex1 & $-2.0769$ & $0.8461,\, (0.7692,0.0000)$ & $1.14\cdot 10^{-8}$ & (3,2,2) & 0.14 \\
	\hline
	MuuQuy2003Ex2 & $0.6389$ & $\begin{array}{l}(0.6111, 0.3889),\\
		(0.0000, 0.00000, 1.8333)\end{array}$ & $7.65\cdot 10^{-8}$ &   (4,2,1) & 0.42 \\
	\hline
	Outrata1990Ex1a & $-8.9172$ & $(1.0316,3.0977),\, (2.5970,1.7929)$ & $9.97\cdot 10^{-7}$ &
	(6,6,3) & 0.40\\
	\hline
	Outrata1990Ex1b & $-7.5785$ & $(0.2788, 0.4748),\, (2.3438, 1.0325)$ & $4.72\cdot 10^{-7}$ &
	(6,6,6) & 0.31\\
	\hline
	Outrata1990Ex1c & $-11.9985$ & $(21.4478, 38.4302),\, (2.9985, 2.9985)$ & $1.90\cdot 10^{-6}$ &
	(6,6,6) & 1.17\\
	\hline
	Outrata1990d & $-3.6000$ & $(2.0000, 0.0000),\, (2.0000, 0.0000)$ & $3.25\cdot 10^{-7}$ &
	(6,6,3) & 0.32\\
	\hline
	Outrata1990e & $-3.9200$ & $(-0.4000, 0.8000),\, (2.0000, 0.0000)$ & $7.91\cdot 10^{-7}$ &
	(6,6,5) & 0.20\\
	\hline
	Outrata19902a & $0.5015$ & $(2.0763),\, (2.9985, 2.9985)$ & $2.86\cdot 10^{-7}$ &
	(6,3,3) & 0.25\\
	\hline
	Outrata19902b & $0.5015$ & $(0.0006),\, (2.9985, 2.9985)$ & $2.30\cdot 10^{-5}$ &
	(6,3,3) & 0.35\\
	\hline
	Outrata19902c & $1.8605$ & $(3.4562),\, (1.7071,2.5685)$ & $2.27\cdot 10^{-6}$ &
	(6,6,6) & 0.56\\
	\hline
	OutrataCervinka2009 & $2.1079\cdot 10^{-9}$ & $(0.0000, 0.2649),\, (-0.1324, -0.1324)$ &
	$8.76\cdot 10^{-6}$ & (3,3,3) & 0.21\\
	\hline
	PaulaviciusEtal2017a & $0.2500$ & $(0.5000),\, (0.0000)$ & $8.33\cdot 10^{-10}$  & (2,2,2) & 0.42 \\
	\hline
	PaulaviciusEtal2017b & $-2.0000$ & $(-1.0000),\, (-1.0000)$ & $2.62\cdot 10^{-8}$ &
	(2,2,1) & 0.09\\
	\hline
	SahinCiric1998Ex2 & $5.0000$ & $(1.0000),\, (3.0000)$ & $3.90\cdot 10^{-8}$ &
	(3,3,3) & 0.17\\
	\hline
	ShimizuAiyoshi1981Ex1 & $100.0000$ & $(10.0000),\, (10.0000)$ & $7.46\cdot 10^{-6}$ &
	(3,2,0) & 0.34\\
	\hline
	ShimizuAiyoshi1981Ex2 & $225.0000$ & $(20.0000, 5.0000),\, (10.0000, 5.0000)$ & $6.65\cdot 10^{-6}$ &
	(4,4,2) & 0.29\\
	\hline
	ShimizuEtal1997a & $16.8919$ & $(0.9459),\, (-0.1622)$ & $3.74\cdot 10^{-9}$ &
	(3,3,2) & 0.10\\
	\hline
	ShimizuEtal1997b & $2250$ & $(11.2500),\,(5.0000)$ & $1.11\cdot 10^{-7}$ &
	(2,2,2) & 0.06\\
	\hline
	SinhaMaloDeb2014TP3 & $-18.6787$ & $(0.0001, 2.0000),\, (1.8750, 0.9062)$ & $8.50\cdot 10^{-8}$ &
	(6, 3, 3) & 0.28\\
	\hline
	SinhaMaloDeb2014TP6 & $-1.2099$ & $(1.8889),\, (0.8889, 0.0000)$ & $6.28\dot 10^{-9}$ &
	(15,9,1) & 0.42\\
	\hline
	SinhaMaloDeb2014TP8 & $0.0000$ & $(0.0000, 30.0000),\, (-10.0000, 10.0000)$ & $3.36\cdot 10^{-6}$&
	(9,9,7) & 0.42\\
	\hline
	TuyEtal2007 & $22.5000$ & $(4.5000),\ (1.5000)$ & $1.34\cdot 10^{-8}$ &
	(3,3,3) & 0.14\\
	\hline
	Vogel2002 & $1.0000$ & $(-2.0000),\, (-2.0000)$ & $5.31\cdot 10^{-8}$ & (1,1,1) & 0.07\\
	\hline
	WanWangLv2011 & $7.5000$ & $\begin{array}{l}(0.5000, 0.5000),\\
		(0.0000, 0.0000, 0.0000)\end{array}$ & $3.90\cdot10^{-9}$ &
	(20,6,6) & 0.10\\
	\hline
	YeZhu2010Ex42 & $1.0000$ & $(1.0000),\, (1.0000)$ & $1.22\cdot 10^{-8}$  & (1,1,1) & 0.10 \\
	\hline
	Yezza1996Ex31 & $0.5000$ & $(3.0000),\, (1.0000)$ & $1.09\cdot 10^{-9}$ &
	(2,2,2) & 0.06\\
	\hline
\end{longtable}
}

\section{Conclusions}\label{sc:con}

In this paper, we propose a PLME-disjunctive decomposition approach for
solving bilevel optimization problems.
For bilevel polynomial optimization which only has
linear lower-level constraints,
we use partial Lagrange multiplier expression to build
the disjunctive decomposition of the KKT system,
and solve it as a set of branch problems.
These branch problems typically have simpler expressions
and are much easier to solve than the original BPOP.
They are equivalent to polynomial optimization problems,
which can be solved efficiently with Moment-SOS relaxations.
Under the existence of feasible extensions,
each branch is either infeasible or has a global optimizer,
for which we give convenient necessary and sufficient conditions for the global optimality,
and sufficient conditions for the local optimality of the original problem.
Our approach performs very efficiently in numerical experiments.

\medskip \noindent
{\bf Acknowledgement}
The first author is partially supported by the NSF grants
DMS-2110780 and DMS-2513254, and the second author is partially supported by NSERC.

\end{document}